%% file: article.tex
\documentclass{preprint}

\journal{Engineering Analysis with Boundary Elements}

\usepackage{amsmath}
\usepackage{amssymb}
\usepackage[]{graphicx}
\usepackage{booktabs}
\usepackage{subfigure}
\usepackage{subfigmat}
\usepackage{color}
\usepackage{url}

\usepackage{symbols-private}

\usepackage[refpage]{nomencl}
\makenomenclature

\begin{document}
\begin{frontmatter}

  \title{A stable and adaptive semi-Lagrangian potential model for
    unsteady and nonlinear ship-wave interactions}

  \author[Sissa]{Andrea Mola}
  \ead{andrea.mola@sissa.it}
  \author[Sissa,cor]{Luca Heltai} 
  \ead{luca.heltai@sissa.it}
  \author[Sissa]{Antonio DeSimone} 
  \ead{desimone@sissa.it}
  \cortext[cor]{Corresponding author, Tel.: +39 040-3787449; Fax: +39 040-3787528}

  \address[Sissa]{SISSA, International School for Advanced Studies, \\ 
    Via Bonomea 265, 34136 Trieste - Italy}

  \begin{abstract}

    We present an innovative numerical discretization of the equations of
    inviscid potential flow for the simulation of three-dimensional,
    unsteady, and nonlinear water waves generated by a ship hull advancing
    in water.  

    The equations of motion are written in a semi-Lagrangian
    framework, and the resulting integro-differential equations are
    discretized in space via an adaptive iso-parametric collocation
    Boundary Element Method, and in time via implicit Backward
    Differentiation Formulas (BDF) with adaptive step size and
    variable order.

    When the velocity of the advancing ship hull is non-negligible, the
    semi-Lagrangian formulation (also known as Arbitrary Lagrangian
    Eulerian formulation, or ALE) of the free-surface equations contains
    dominant transport terms which are stabilized with a Streamwise Upwind
    Petrov-Galerkin (SUPG) method.

    The SUPG stabilization allows automatic and robust adaptation of
    the spatial discretization with unstructured quadrilateral
    grids. Preliminary results are presented where we compare our
    numerical model with experimental results on a Wigley hull
    advancing in calm water with fixed sink and trim.
  \end{abstract}
\end{frontmatter}

\section{Introduction}

The use of computational tools to predict hydrodynamic performances of
ships has gained a lot of popularity in recent years. Models based on
potential flow theory have historically been among the most successful
to assess the wave pattern around a ship hull in presence of a forward
ship motion.

In this framework, the assumptions of irrotational flow and inviscid
fluid reduce the Navier Stokes incompressibility constraint and
momentum balance equations to the Laplace's and Bernoulli's equations,
defined on a moving and \textit{a-priori} unknown domain.

This boundary value problem is tackled with a Mixed
Eulerian--Lagrangian approach, in which the Eulerian field equations
are solved to obtain the fluid velocities, which are then used to
displace in a Lagrangian way the free-surface, and update the
corresponding potential field values \cite{tanizawa2000}. In this
framework, the Eulerian problem is expressed in boundary integral
form, and it is typically discretized using the Boundary Element
Method (BEM).
%
The velocity field and Bernoulli's equation provide a kinematic
boundary condition for the Lagrangian evolution of the free-surface,
and a dynamic boundary condition for the evolution of the potential
field.

Numerical treatment of the Lagrangian step usually relies on accurate
reconstructions of the position vector and of potential field
gradients. In presence of a forward motion, these reconstructions may
lead to instability in the time advancing scheme for the free-surface
discretization, as well as for the corresponding potential field
values. A smoothing technique is often adopted to reduce the
sensitivity of the discretization on the reconstruction of the full
velocity field, at the cost of introducing an artificial viscosity in
the system.

An alternative cure was presented by Grilli et al. \cite{grilli2001}
who developed a high order iso-parametric BEM discretization of a
Numerical Wave Tank to simulate overturning waves up to the breaking
point on arbitrarily shaped bottoms. The use of a high order
discretization bypasses the problem of reconstructing the gradients,
and is very reliable when the numerical evolution of the free-surface
is done in a purely Lagrangian way.

Ship hydrodynamics simulations, however, are typically carried out in
a frame of reference attached to the boat, requiring the presence of a
water current in the simulations. In this case a fully Lagrangian
approach leads to downstream transportation of the free-surface
nodes, or to their clustering around stagnation points, ultimately
resulting in blowup of the simulations (see, for example,
\cite{scorpioPhD}).

Beck~\cite{beck1994} proposed alternative \emph{semi-Lagrangian} free
surface boundary conditions, under the assumption that the surface
elevation function is single-valued. Employing such conditions, it is
possible to prescribe arbitrarily the horizontal velocity of the free
surface nodes, while preserving the physical shape of the free-surface
itself.

However, if there is a significant difference between the water
current speed and the horizontal nodes speed, stability issues may
arise \cite{scorpioPhD}, \cite{RyuKimLynett2003}. For this reason, the
semi-Lagrangian free-surface boundary conditions proposed by Beck have
been in most cases employed imposing nodes longitudinal speeds equal
to the water current ones \cite{scorpioPhD}, thus reducing the
instability at the cost of more complex algorithms for mesh
management.

More recently, Sung and Grilli \cite{SungGrilli2008} applied an
alternative method, combining semi-Lagrangian and Lagrangian free
surface boundary conditions to the problem of a pressure perturbation
moving on the water surface
. The semi-Lagrangian free-surface boundary conditions were used
also in the work of Kjellberg, Contento and Jansson
\cite{KjelContJans-2011}, where free-surface instabilities are
avoided by choosing an earth-fixed reference frame, in which no
current speed is needed. The drawback of this choice is that in such
frame the ship moves with a specified horizontal speed, and the
computational grid needs to be constantly regenerated to cover the
region surrounding the hull with an adequate number of cells.

The purpose of this work is to present a new stabilized
semi-Lagrangian potential model for the simulation of three
dimensional unsteady nonlinear water waves generated by a ship hull
advancing in calm water. The resulting integro-differential boundary
value problem is discretized to a system of nonlinear
differential-algebraic equations, in which the unknowns are the
positions of the nodes of the computational grid, along with the
corresponding values of the potential and of its normal
derivative.

Time advancing of the nonlinear differential-algebraic system is
performed using implicit Backward Differentiation Formulas (BDF) with
variable step size and variable order, 
 implemented in the
framework of the open source library SUNDIALS \cite{sundials2005}
. The
collocated and iso-parametric BEM discretization of the Laplace's
equation has been implemented using the open source C++ library
deal.II \cite{BangerthHartmannKanschat2007}.

The computational grid, composed by quadrilateral cells of arbitrary
order, is adapted in a geometrically consistent way
(see~\cite{BonitoNochettoPauletti-2010-a}) via an
\textit{a-posteriori} error estimator based on the jump of the
solution gradient along the cell internal boundaries. Even when low
order boundary elements are used, accurate estimations of the position
vector and potential gradients on the free-surface are recovered by a
Streamwise Upwind Petrov--Galerkin (SUPG) projection, which is used to
stabilize the transport dominated terms. The SUPG projection is
strongly consistent, and does not introduce numerical dissipation in
the equations. This allows the use of robust unstructured grids, which
can be generated and managed on arbitrary hull geometries in a
relatively simple way.

The test case considered in this paper is that of a Wigley hull
advancing at constant speed in calm water. The simulations have been
performed using bilinear elements. The arbitrary horizontal velocity
specified in the semi-Lagrangian free-surface boundary condition is
chosen so the nodes have  null longitudinal  velocity with respect
to the hull. The numerical results obtained imposing six different
Froude numbers are finally compared with experimental results reported
in \cite{mccarthy1985}, to assess the accuracy of the model.

\section{Three-dimensional potential model}

The equations of motion that describe the velocity and pressure
fields $\vb$ and $p$ of a fluid region around a ship hull are the
incompressible Navier--Stokes equations. Such equations, written
in the moving domain $\Omega(t) \subset\mathbb{R}^3$ (a simply
connected region of water surrounding the ship hull) read:


\begin{subequations}
  \label{eq:incompressible-euler}
  \begin{alignat}{3}
    \label{eq:conservation-momentum}
    & \rho \left( \frac{\partial \vb}{\partial t} + \vb \cdot \nablab
    \vb \right)= \nablab\cdot\sigmab + \bb \qquad & \text{ in } \Omega(t)\\
    \label{eq:incompressibility}
    & \nablab \cdot \vb = 0 & \text{ in } \Omega(t)\\
    \label{eq:freesurface-condition}
    & \sigmab\cdot\nb = p_a \nb & \text{ on } \Gamma^w(t)\\
    \label{eq:non-penetration-condition}
    & \vb = \vb_g & \text{ on } \Gamma(t)\setminus
    \Gamma^w(t)
  \end{alignat}
\end{subequations}
where $\rho$ is the (constant) fluid density, $\bb$ are external body
forces (typically gravity and possibly other inertial forces due to a
non inertial movement of the reference frame), $\sigmab = -p\Ib +
\mu(\nablab\vb+\nablab\vb^T)$ is the stress tensor for an
incompressible Newtonian fluid, $\Gamma(t) := \partial\Omega(t)$ is
the boundary of the region of interest, and $\nb$ is the outer normal
to the boundary $\Gamma(t)$. On the free-surface $\Gamma^w(t)$, the
air is assumed to exert a constant atmospheric pressure $p_a$ on the
underlying water, and we neglect shear stresses due to the wind. On
the other boundaries of the domain, the prescribed velocity $\vb_g$ is
either equal to the ship hull velocity, or to a given velocity field
(for inflow and outflow boundary conditions far away from the ship
hull itself).

\nomenclature{$t$}{Time}
\nomenclature{$\vb$}{Fluid velocity}
\nomenclature{$\vb_g$}{Velocity (given) boundary conditions}
\nomenclature{$\sigmab$}{Fluid Cauchy stress}
\nomenclature{$\Omega(t)$}{Fluid domain}
\nomenclature{$\Gamma(t)$}{Boundary of the fluid domain}
\nomenclature{$\Gamma^w(t)$}{free-surface}
\nomenclature{$p$}{Pressure}
\nomenclature{$p_a$}{Atmospheric pressure}
\nomenclature{$\rho$}{Fluid density}
\nomenclature{$\bb$}{Body forces}

Equation~(\ref{eq:conservation-momentum}) is usually referred to as
the momentum balance equation, while (\ref{eq:incompressibility}) is
referred to as the incompressibility constraint, or continuity
equation.

In the flow field past a slender ship hull, vorticity is confined to
the boundary layer region and to a thin wake following the boat: in
these conditions, the assumption of irrotational and non viscous flow
is fairly accurate. The neglected viscous effects can be recovered, if
necessary, by other means such as empirical algebraic formulas,
or ---better--- by the interface with a suitable boundary layer model.

For an irrotational and invishid flow, the velocity field $\vb$
admits a representation through a scalar potential function
$\Phi(\xb,t)$, namely
\begin{equation}
\label{eq:potential-definition}
\vb = \nablab\Phi \qquad\qquad    \text{in} \ \ \Omega(t).
\end{equation}
\nomenclature{$\Phi$}{Total velocity potential}
In this case, the equations of motion simplify to the unsteady
Bernoulli equation and to the Laplace equation for the flow potential:

\begin{subequations}
  \label{eq:incompressible-euler-potential}
  \begin{alignat}{3}
    \label{eq:bernoulli}
    & \frac{\Pa\Phi}{\Pa t}+\frac{1}{2}\left|\nablab\Phi\right|^2
      +\frac{p-p_a}{\rho}+\bbeta
      \cdot \xb = C(t) \qquad & \text{ in } \Omega(t)\\
    \label{eq:incompressibility-potential}
    & \Delta \Phi = 0 & \text{ in } \Omega(t)
  \end{alignat}
\end{subequations}
\nomenclature{$\bbeta$}{Body force potential}
\nomenclature{$\xb$}{Vector coordinates in real space}
where $C(t)$ is an arbitrary function of time, and we have assumed
that all body forces could be expressed as $\bb = -\nablab (\bbeta
\cdot \xb)$, i.e., they are all of potential type. This is true for
gravitational body forces and for inertial body forces due to uniform
accelerations along fixed directions of the frame of reference.

In this framework, the unknowns of the problem $\Phi$ and $p$ are
uncoupled, and it is possible to recover the pressure by postprocessing
the solution of the Poisson
problem~(\ref{eq:incompressibility-potential}) via Bernoulli's
Equation~(\ref{eq:bernoulli}). 

Since arbitrary uniform variations of the potential field produce the
same velocity field (i.e., $\nablab(\Phi(\xb, t)+C(t)) =
\nablab\Phi(\xb, t)$), we can set $C(t)=0$, and we can combine
Bernoulli equation~(\ref{eq:bernoulli}) and the dynamic boundary
condition on the free-surface~(\ref{eq:freesurface-condition}) to
obtain an evolution equation for the potential field $\Phi(\xb,t)$ on
the free-surface $\Gamma^w(t)$.

The shape of the water domain $\Omega(t)$ is time-dependent and it is
part of the unknowns of the problem. The free-surface $\Gamma^w(t)$
should move following the fluid velocity $\vb$.  In ship
hydrodynamics, however, it is desirable to maintain the frame of
reference attached to the boat, and study the problem in a domain
which is neither fixed nor a material subdomain. Indeed, its evolution
is not governed solely by the fluid motion, but also by the motion of
the reference frame and by the motion of the stream of
water. In addition, it has to comply to the free-surface boundary
$\Gamma^w(t)$ which is the result of the dynamic boundary
condition~(\ref{eq:freesurface-condition}).

A possible solution is to introduce an intermediate frame of
reference, called Arbitrary Lagrangian Eulerian (ALE) (see, for
example, \cite{formaggia2009cardiovascular} and~\ref{sec:ale}). In the
context of potential free-surface flows, this approach is also known
as the \emph{semi-Lagrangian}
approach~\cite{BeckCaoLee-1993-a,scorpioPhD}.

\subsection{Perturbation potential and boundary conditions}

\begin{figure}[htb!]        
\centerline{
  \ifpdf
  \resizebox{0.8\textwidth}{!}{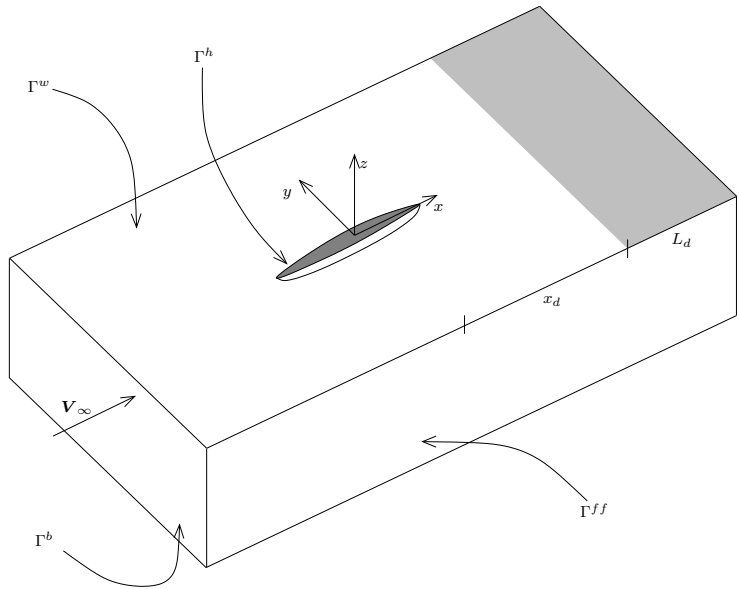}
  \else
  \resizebox{0.8\textwidth}{!}{\input{figures/dominio_barca.pstex_t}}
  \fi
}
\caption{The domain $\Omega(t)$ and the different regions
         in which its boundary $\Pa\Omega(t)$ is split. The grey area
         behind the hull $\Gamma^{h}(t)$, is the absorbing beach portion
         of the free-surface region $\Gamma^{w}(t)$.  
         \label{fig:wholeDomain}}	
\end{figure} 

We solve Problem (\ref{eq:incompressible-euler-potential}) on the
region $\Omega(t)$ around the boat, and we move the local frame of
reference with horizontal velocity $\Vb_{\sys}(t)$ which coincides
with the horizontal boat velocity (which we assume as known). An
additional horizontal water stream velocity $\Vb_{\stream}(t)$,
expressed in a global (earth-fixed) reference frame is added to the
problem.

Uniform accelerations of the reference frame $\ab_{\sys}$ can be taken into account
by incorporating them in the body force term
\begin{equation}
  \label{eq:forcing-term}
  \bbeta\cdot\xb = \ab_{\sys}\cdot\xb + gz,
\end{equation}
and it is convenient to split the potential $\Phi$ into the sum between
a mean flow potential (due to the stream velocity and to the frame of
reference velocity) and the so called \emph{perturbation potential}
$\phi$ due to the presence of the ship hull, namely

\nomenclature{$g$}{Gravity acceleration}
\nomenclature{$\ab_{\sys}$}{Acceleration of the reference frame}

\begin{subequations}
  \begin{alignat}{1}
    \Phi(\xb,t) &= \big(\Vb_{\stream}(t)-\Vb_{\sys}(t)\big)\cdot\xb +
    \phi(\xb,t)\\
    \vb(\xb,t) = \nablab\Phi(\xb,t) &= \Vb_{\stream}(t)-\Vb_{\sys}(t) +
    \nablab\phi(\xb,t)\\
    \frac{\Pa\Phi}{\Pa t}(\xb,t) &=
    \big(\ab_{\stream}(t)-\ab_{\sys}(t)\big)\cdot\xb +
    \frac{\Pa\phi}{\Pa t}(\xb,t).
  \end{alignat}
  \label{eq:potSplit}
\end{subequations}
\nomenclature{$\phi$}{Perturbation potential}
\nomenclature{$\Vb_{\sys}$}{Translation speed of the frame of
  reference}
\nomenclature{$\Vb_{\stream}$}{Water stream velocity}
The perturbation potential still satisfies Poisson equation

\begin{equation}
\Delta\phi = 0 \qquad\qquad \text{in} \ \ \Omega(t),
\end{equation}
and in practice it is convenient to solve for $\phi$, and obtain the
total potential $\Phi$ from equation (\ref{eq:potSplit}).

In Figure~\ref{fig:wholeDomain} we present a sketch of the domain
$\Omega(t)$, with the explicit splitting of the various parts of the
boundary $\Gamma(t)$. On the boat hull surface, the non-penetration
condition takes the form

\begin{equation}
\phin := \nablab\phi\cdot \nb= \nb\cdot(\Vb_{\boat}-\Vb_{\stream})
\qquad\qquad \text{on} \ \ \Gamma^{h}(t),
\end{equation}
when expressed in terms of the perturbation potential and
$\Vb_{\boat}$ is the (given) boat velocity. 

\nomenclature{$\Gamma^{h}(t)$}{Hull surface}

On the ---horizontal--- bottom of the water basin $\Gamma^{b}$, the
non penetration condition is also applied, namely

\begin{equation}
\phin = 0 \qquad\qquad \text{on} \ \ \Gamma^{b}(t).
\end{equation}
A possible condition for the ---vertical--- far field boundary is the
homogeneous Neumann condition

\nomenclature{$\Gamma^{b}(t)$}{Bottom surface of the basin}

\begin{equation}
\phin = 0 \qquad\qquad \text{on} \ \ \Gamma^{ff}(t).
\end{equation}
The most remarkable limit of such condition is the fact that it
reflects wave energy back in the computational domain, thus limiting
the simulation time. Different boundary conditions can be used to let
the wave energy flow outside the domain. In~\ref{sec:numerical-beach},
we explain in detail the procedure used in our computations.

\nomenclature{$\Gamma^{ff}(t)$}{Far field surface of the basin}

Applying the potential splitting~(\ref{eq:potSplit}) to
Problem~(\ref{eq:incompressible-euler-potential}), and summarizing all
boundary conditions, we obtain the perturbation potential formulation

\begin{subequations}
  \label{eq:incompressible-euler-perturbation-potential}
  \begin{alignat}{2}
    \label{eq:incompressibility-potential-perturbation-phi}
    & \Delta \phi = 0 & \text{ in } \Omega(t)\\
    \label{eq:fsDynamicLag}
    & \frac{\Pa\phi}{\Pa t} +\vb\cdot\nablab\phi= -gz +
    \ab_{\stream}\cdot\xb+ \frac{1}{2}|\nablab\phi|^2 & -
    \mu\phin \qquad \text{ on } \Gamma^w(t) \\
    &\phi(\xb, 0) = \phi_0(\xb) &\text{ on } \Gamma^w(0)\label{eq:initial-condition-phi} \\
    \label{eq:non-penetration-condition-hull}
    &\phin =
    \nb\cdot(\Vb_{\boat}-\Vb_{\stream}) & \text{ on }
    \Gamma^{h}(t)\\
    &\phin = 0 & \text{ on }
    \Gamma^b(t)\cup\Gamma^{ff}(t)\\
    \label{eq:fsKinematicLag}
    &\frac{D \Gamma^w(t)}{D t} = \vb &
    \text{ on }\Gamma^w(t),
  \end{alignat}
\end{subequations}
where $\vb = \Vb_{\stream} - \Vb_{\sys}+ \nablab\phi$.
Equation~(\ref{eq:fsKinematicLag}) is a kinematic boundary condition
and indicates that the free-surfaces follows the velocity of the
fluid, while Equation~(\ref{eq:fsDynamicLag}) is a dynamic Dirichlet
boundary condition for the Poisson
problem~(\ref{eq:incompressibility-potential-perturbation-phi}) at
each time $t$, whose initial condition is given by
Equation~(\ref{eq:initial-condition-phi}). A detailed explanation of
the role of the term $\mu \phin$ in Equation~\eqref{eq:fsDynamicLag}
is given in~\ref{sec:numerical-beach}.

\nomenclature{$\Vb_{\boat}$}{Velocity of the boat}




Introducing a non-physical reference domain (the Arbitrary Lagrangian
Eulerian domain $\tilde \Omega$) and an arbitrary deformation
$\pale(\xale,t)$ such that $\pale(\tilde \Omega, t) = \Omega(t)$ as
explained in~\ref{sec:ale}, Equations~\eqref{eq:fsDynamicLag}
and~\eqref{eq:fsKinematicLag} can be rewritten in terms of ALE
derivatives $\delta/\delta t := \Pa / \Pa t + \wb\cdot \nablab$ as

\begin{subequations}
  \begin{alignat}{2}
    \label{eq:fsDynamicSemiLag}
    & \frac{\delta\phi}{\delta t} +(\vb-\wb)\cdot\nablab\phi= -g\eta
    + \ab_{\stream}\cdot\xb+ \frac{1}{2}|\nablab\phi|^2 & -
    \mu\phin \qquad \text{ on } \Gamma^w(t) \\
    \label{eq:condition-on-ALE-velocity}
    &\wb\cdot \nb = \vb\cdot \nb &\text{ on }
    \Gamma(t),
  \end{alignat}
\end{subequations}
where $\wb(\pale(\xale,t), t) =\Pa \pale(\xale,t)/\Pa t$ is the
velocity of the domain deformation expressed in Eulerian
coordinates (see~\ref{sec:ale} for the details). 

We remark that, whenever $\wb=\vb$, we recover a fully Lagrangian
formulation: the ALE motion is, in this case, following the particles
motion and $\delta/\delta t \equiv D/Dt$. Similarly, if the domain is
fixed, and we set $\wb=0$, we recover the classical Eulerian
formulation~(\ref{eq:incompressible-euler-perturbation-potential}) on
fixed domains, and $\delta/\delta t \equiv \Pa/\Pa t$.

If the free-surface can be seen as the graph of a single valued
function $\eta(x,y,t)$ of the horizontal components $x$ and $y$, then
\begin{equation}
  \label{eq:definition-eta}
  \xb = (x,y, \eta(x,y,t)) \qquad \text{ on }\Gamma^w(t),
\end{equation}
and, for a material particle on the free-surface, we have
\begin{equation}
  \label{eq:particle-eta}
  \plag(\xlag,t)\cdot\eb_z = \eta(\plag(\xlag, t),t),
\end{equation}
where we define $\eta(\xb,t) :=
\eta(\xb\cdot\eb_x,\xb\cdot\eb_y,t)$. Taking the time derivative of
Equation~(\ref{eq:particle-eta}) we get
\begin{equation}
  \frac{\Pa \plag}{\Pa t}(\xlag,t)\cdot\eb_z = \frac{\Pa \eta}{\Pa
    t}(\plag(\xlag, t),t) + \frac{\Pa \plag }{\Pa t}(\xlag, t)\cdot
  \nablab \eta(\plag(\xlag, t),t)  \qquad \text{ on } \hat\Gamma^w,
\end{equation}
\nomenclature{$\eta(x,y,t)$}{Elevation of the free-surface}
that is, in Eulerian form:
\begin{equation}
  \label{eq:particle-veolocity-on-freesurface-eulerian}
  \vb_z = \frac{\Pa \eta }{\Pa
    t} +\vb\cdot \nablab \eta = \frac{D \eta}{D t}\qquad \text{ on } \Gamma^w(t),
\end{equation}
where $\Pa\eta/\Pa z \equiv 0$. Proceeding in the same way for the ALE
deformation and the ALE velocity of the domain, we get
\begin{equation}
  \label{eq:ale-veolocity-on-freesurface-eulerian}
  \wb_z = \frac{\Pa \eta}{\Pa
    t} +\wb\cdot \nablab \eta = \frac{\delta \eta}{\delta t}\qquad \text{ on } \Gamma^w(t).
\end{equation}
Isolating $\Pa \eta/\Pa t$ in
Equation~(\ref{eq:particle-veolocity-on-freesurface-eulerian}) and
substituting in
Equation~(\ref{eq:ale-veolocity-on-freesurface-eulerian}), we get an
alternative expression for
Condition~(\ref{eq:condition-on-ALE-velocity})\footnote{%
  Conditions~(\ref{eq:condition-on-ALE-velocity})
  and~(\ref{eq:fsKinematicSemiLag}) are equivalent in this
  framework. This comes from the observation that, when $\Gamma^w(t)$
  is the graph of the function $\eta(x,y,t)$, then the normal to the
  surface $\Gamma^w(t)$ itself is given by a vector proportional to
  $\eb_z-\nablab\eta$. Substituting in
  Equation~(\ref{eq:condition-on-ALE-velocity}) gives immediately
  Equation~(\ref{eq:fsKinematicSemiLag}).},
\begin{equation}
  \frac{\delta\eta}{\delta t} +(\vb-\wb)\cdot\nablab\eta =
    \vb\cdot \eb_z
    \qquad \text{ on }\Gamma^w(t),\label{eq:fsKinematicSemiLag}
\end{equation}
which is valid for nonbreaking waves in which $\eta(x,y,t)$ is
single-valued. 

Equation~(\ref{eq:fsKinematicSemiLag}) is the kinematic boundary
condition for the evolution of the unknown free-surface elevation
$\eta(x,y,t)$ which is often found in the literature of
semi-Lagrangian methods for potential free-surface
flows~\cite{beck1994,BeckCaoLee-1993-a,scorpioPhD}.

Equation (\ref{eq:fsKinematicSemiLag}) is rather general, and is valid
for \emph{arbitrary} values of horizontal ALE velocities. Suitable
values of $\Vb_{\stream}$ and $\Vb_{\sys}$ can be plugged into the
velocity field $\vb = (\Vb_\stream-\Vb_\sys +\nablab\phi)$ to specify
them for the desired reference frame. For instance, setting
$\Vb_{\stream} = \Vb_{\sys} = 0$, and
$\wb=\left(0,0,\frac{\delta\eta}{\delta t}\right)$, one obtains

\begin{eqnarray}
\frac{\delta\eta}{\delta t} &=& \frac{\Pa\phi}{\Pa z}+\nablab\eta\cdot\left(\wb -\nablab\phi\right)
                               = \frac{\Pa\phi}{\Pa z}-\frac{\Pa\phi}{\Pa x}\frac{\Pa\eta}{\Pa x}
                                 -\frac{\Pa\phi}{\Pa y}\frac{\Pa\eta}{\Pa y}
\label{eq:fsKinematicContento}\\
\frac{\delta\phi}{\delta t} &=& -g\eta +
\frac{1}{2}|\nablab\phi|^2 + \nablab\phi\cdot\left(\wb -\nablab\phi\right) 
                               =  -g\eta + \frac{\Pa\phi}{\Pa z}\frac{\delta\eta}{\delta t}-\frac{1}{2}|\nablab\phi|^2
\label{eq:fsDynamicContento},
\end{eqnarray}  
which are the semi-Lagrangian equations written in an earth fixed reference frame, and null stream velocity
used in \cite{KjelContJans-2011}.

In this work, we choose instead to solve the problem in a coordinate
system attached to the ship hull. We move the reference frame
according to the horizontal average velocity of the boat, that is we
set $\Vb_\sys = (\overline\Vb_\boat\cdot\eb_x)\eb_x =: -\Vb_\infty =
\left(-V_\infty,0,0\right)$ and we assume $\Vb_{\stream} = 0$. With
these values, Equations (\ref{eq:fsKinematicSemiLag}) and
(\ref{eq:fsDynamicSemiLag}) take the form

\begin{eqnarray}
\frac{\delta\eta}{\delta t} &=& \frac{\Pa\phi}{\Pa z}+\nablab\eta\cdot\left(\wb -\nablab\phi-\Vb_\infty\right)
\label{eq:fsKinematicBeck}\\
\frac{\delta\phi}{\delta t} &=& -g\eta +
\frac{1}{2}|\nablab\phi|^2 + \nablab\phi\cdot\left(\wb -\nablab\phi-\Vb_\infty\right)-\mu\phin,
\label{eq:fsDynamicBeck}
\end{eqnarray}  
which coincide with the ones proposed by Beck \emph{et al.}
\cite{BeckCaoLee-1993-a}, and in which the points on the free-surface
move with an a priori arbitrary horizontal speed in the boat
reference frame.

\nomenclature{$\Vb_{\infty}$}{Horizontal average velocity of the boat}

\subsection{Boundary integral formulation}

While
Equation~(\ref{eq:incompressibility-potential-perturbation-phi})
is time-dependent and defined in the entire domain $\Omega(t)$, we are
really only interested in its solution on the boundary $\Gamma(t)$. 

Knowledge of the solution on the free-surface part of the boundary is
used to advance the shape of the domain in time, while Bernoulli's
equation~(\ref{eq:bernoulli}) is used to recover the pressure
distribution on the ship hull by postprocessing the solution on
$\Gamma^h(t)$.

In other words, at any given time instant $\overline{t}$ we want to
compute $\phi$ satisfying
\begin{subequations}
\label{eq:laplace-problem}
\begin{alignat}{2}
-\Delta\phi &= 0 & \text{in} \ \ \Omega(\overline{t}) \\
\phi &= \overline{\phi} & \text{on} \ \ \Gamma^{w}(\overline{t}) \\
\phin &= \overline\phin   & \qquad\text{on} \ \ \Gamma^{h}(\overline{t})\cup\Gamma^{b}(\overline{t})\cup \Gamma^{ff}(\overline{t})  
\end{alignat}
\end{subequations}
where $\overline{\phi}$ is the potential on the free-surface at the
time $\overline{t}$, and $\overline\phin$ is
equal to zero on $\Gamma^b(\overline{t} )\cup\Gamma^{ff}(\overline{t}
)$ and to $\big(\Vb_\boat(\overline{t})-\Vb_\infty(\overline{t})\big)\cdot\nb$ on
$\Gamma^{h}(\overline{t} )$.

This is a purely spatial boundary value problem, in which time appears
only through boundary conditions and through the shape of the time
dependent domain.

The solution of this boundary value problem allows for the computation of
the full potential gradient on the boundary $\Gamma(t)$, which is what
is required in the dynamic and kinematic boundary conditions to
advance the time evolution of both $\phi$ and $\eta$.

Using the second Green identity
\begin{equation}
  \int_{\Omega}
  (-\Delta u)v\,dx + \int_{\partial\Omega} \frac{\partial u}{\partial
    n}v \d s
  = 
  \int_{\Omega}
  (-\Delta v)u\,dx + \int_{\partial\Omega} \frac{\partial v}{\partial
    n} u \d s,
\end{equation}
a solution to system (\ref{eq:laplace-problem}) can be expressed in terms
of a boundary integral representation only, via convolutions with
fundamental solutions of the Laplace problem.

We call $G$ the \emph{fundamental solution}, i.e., the function
\begin{equation}
G(\rb)  =  \frac{1}{4\pi |\rb|},
\end{equation}
\nomenclature{$G$}{Fundamental solution}
which is the distributional solution of
\begin{equation} \label{eq:poissonGreen}
  \begin{split}
    &-\Delta G(\xb-\xb_0) = \delta(\xb_0) \qquad\qquad
     \text{ in } \Re^3\\
    &\lim_{|\xb|\to\infty} G(\xb-\xb_0) = 0,
  \end{split}
\end{equation} 
\nomenclature{$\delta(\xb_0)$}{Dirac delta distribution centered in $\xb_0$}
where $\delta(\xb_0)$ is the Dirac delta distribution centered in $\xb_0$.

If we select $\xb_0$ to be inside $\Omega(t)$, use the defining property
of the Dirac delta and the second Green identity, we obtain
\begin{multline}
\phi(\xb_0,t) =
\int_{\Omega(t)}\left[-\left(\Delta G(\xb-\xb_0)\right) \phi(\xb,t)\right] \d\Omega  = \\
\int_{\Gamma(t)} \left[(\nablab\phi(\xb,t)\cdot\nb)G(\xb-\xb_0) -
(\nablab G(\xb-\xb_0)\cdot\nb)\phi(\xb,t)\right]\d\Gamma.
\end{multline}

In the limit for $\xb_0$ touching the boundary $\Gamma(t)$, the
integral on the right hand side will have a singular argument, and
should be evualated according to the Cauchy principal value. This
process yields the so called \emph{Boundary Integral Equation} (BIE)
\begin{equation}
\label{eq:bie}
\alpha\phi  = 
\int_{\Gamma(t)} \left[\phin G -
\frac{\Pa G}{\Pa n}\phi\right]\d\Gamma \qquad
  \text{ on } \Gamma(t),
\end{equation}
where $\alpha(\xb,t)$ is the fraction of solid angle $4\pi$ with which
the domain $\Omega(t)$ is seen from $\xb$ and the gradient of the
fundamental solution is given by
\begin{equation}
\nablab G(\rb)\cdot \nb  = 
-\frac{\rb\cdot\nb}{4\pi \left|\rb\right|^3}.
\end{equation} 
\nomenclature{$\alpha(\xb)$}{Fraction of solid angle with which
  $\Omega$ is seen from $\xb$}

The function $\alpha(\xb,t)$ can be computed by noting that the
constant function $1$ is a solution to the Laplace equation with zero
normal derivative, and therefore it must be
\begin{equation}
  \label{eq:computation-alpha}
  \alpha = -\int_{\Gamma(t)} \frac{\Pa G}{\Pa n} d\Gamma\qquad
  \text{ on } \Gamma(t),
\end{equation}
in the Cauchy principal value sense.

With Equation~(\ref{eq:bie}), the continuity equation has been
reformulated as a boundary integral equation of mixed type defined on
the moving boundary $\Gamma(t)$, where the main ingredients are the
perturbation potential $\phi(\xb,t)$ and its normal derivative
$\displaystyle\phin(\xb,t)$.

The domain deformation $\pb(\xb,t)$ on the free-surface takes the form
\begin{equation}
  \label{eq:def-ale-motion-gamma}
  \pb(\xb,t) = (x,y,\eta(x,y,t)) \qquad \text{ on } \Gamma^w(t), 
\end{equation}
and one has to solve an additional boundary value problem to uniquely
determine the full ALE motion $\pale(\xale,t)$.

\subsection{Arbitrary Lagrangian Eulerian motion}
\label{sec:arbitr-lagr-euler}

When the arbitrary Lagrangian Eulerian formulation is used in the
finite element framework (see, for
example,~\cite{formaggia2009cardiovascular}), the restriction to the
boundary $\Gamma(t)$ of the deformation $\pale(\xale,t)$ is either
known, or entirely determined by the equations of motion. In this
case, an additional boundary value problem needs to be solved to
recover the domain deformation in the interior of $\Omega(t)$ starting
from the Dirichlet values on the boundary.

Our situation is slightly different, since only the \emph{normal}
component of the motion is given on the boundary $\Gamma(t)$, and we
are not really interested in finding a domain motion in the interior
of $\Omega(t)$.

In the dynamic and kinematic boundary conditions
(\ref{eq:fsKinematicBeck}) and~(\ref{eq:fsDynamicBeck}) we have the
freedom to choose an ALE motion arbitrarily, as long as the shape of
$\Gamma(t)$ is preserved. In analogy to what is done in the finite
element framework, we construct an additional boundary value problem
to determine uniquely $\pale(\xale,t)$.

A typical choice in the finite element framework is based on linear
elasticity theory, and requires the solution of an additional Laplace
problem on the coordinates $\pale(\xale,t)$, or, in some cases, a
bi-Laplacian. This procedure can be generalized to surfaces embedded
in three-dimensional space via the Laplace-Beltrami operator. 

The Laplace-Beltrami operator can be constructed from the surface
gradient $\nablab_s g(\xb,t)$, defined as
\begin{equation}
  \label{eq:surface-gradient}
    \nablab_s a(\xb,t)  := \nablab
    \overline{a}-(\nablab\overline{a}\cdot\nb)\nb, \qquad
    \forall \overline a \text{ s.t. } \overline a = a(\xb,t) \quad\text{ on }\Gamma(t),
\end{equation}
\nomenclature{$\nablab_s$}{Surface gradient}
where $\overline{a}$ is an arbitrary smooth extension of $a(\xb,t)$ on
a tubular neighborhood of
$\Gamma(t)$. Definition~(\ref{eq:surface-gradient}) is independent on
the extension used (see, for example,
\cite{delfour2010shapes}). Similarly, we indicate with
$\tilde\nablab_s$ the surface gradient computed in the reference
domain $\tilde\Gamma$, with the same definition as
in~(\ref{eq:surface-gradient}), but replacing $\xb$ with $\xale$, and
performing all differential operators in terms of the independent
variable $\xale$ instead of $\xb$.

If we indicate with $\nablab_s\cdot$ the surface divergence (i.e., the
trace of the surface gradient $\nablab_s$), then the surface Laplacian
$\Delta_s$ and $\tilde\Delta_s$ on $\Gamma(t)$ and on $\tilde\Gamma$,
are given by $\Delta_s:= \nablab_s\cdot\nablab_s$, and by $\tilde\Delta_s:=
\tilde\nablab_s\cdot\tilde\nablab_s$.

We use the shorthand notation $\gamma^{a,b}(t)$
to indicate the intersection between $\Gamma^a(t)$ and $\Gamma^b(t)$,
that is,
\begin{equation}
  \gamma^{a,b}(t) = \overline{\Gamma^a(t)}\cap \overline{\Gamma^b(t)}
  \qquad a \neq b,
\end{equation}
where $a,b$ are either $w$, $h$, $b$ or $ff$. We indicate with
$\gamma(t)$ the union of all curves $\gamma^{a,b}(t)$. 

\nomenclature{$\gamma^{a,b}$}{Curve intersection between surfaces
  $\Gamma^a$ and $\Gamma^b$}

The curve $\gamma^{w,h}(t)$ is usually referred to as the
\emph{waterline} on the hull of the ship. On $\gamma^{w,h}(t)$, the
domain velocity $\wb$ has to satisfy the kinematic boundary condition
for \emph{both} the free-surface and the ship hull:
\begin{equation}
  \label{eq:velocity-on-waterline}
  \begin{aligned}
    \wb\cdot\nb^w &= \vb\cdot\nb^w \qquad& \text{ on } \gamma^{w,h}(t)&\\
    \wb\cdot\nb^h &= 0 & \text{ on } \gamma^{w,h}(t)&,
  \end{aligned}
\end{equation}
where $\nb^w$ is the normal to the free-surface and $\nb^h$ is the
normal to the hull surface. When both conditions are enforced, $\wb$
is still allowed to be arbitrary along the direction tangent to the
waterline. 

There are several options to select the tangent velocity $\wb_t$
defined as
\begin{equation}
  \wb_t := \wb\cdot(\nb^h\times\nb^w) = \wb\cdot \tb.
\end{equation}
A natural possibility is to choose zero tangential velocity. Other
choices are certainly possible, and may be preferable, for example, if
one would like to cluster computational nodes in regions where the
curvature of the waterline is higher. In the experiments we present,
the tangential velocity is always set to zero.

Conditions~(\ref{eq:velocity-on-waterline}) and zero tangential
velocity, uniquely determine an evolution equation for the ALE
deformation on $\tilde\gamma^{w,h}$. Here we summarize all boundary
conditions for the evolution of $\pale(\xale,t)$ on the entire
$\tilde\gamma$:
\begin{subequations}
  \label{eq:w-cond-gamma}
  \begin{align}
    \wb\cdot\nb^w &= \vb\cdot\nb^w \qquad& \text{ on } \gamma^{w,h}(t)&\nonumber\\
    \wb\cdot\nb^h &= 0 & \text{ on } \gamma^{w,h}(t)& \label{eq:w-cond-gamma-w-h}\\
    \wb\cdot(\nb^h\times\nb^w) &= 0 & \text{ on } \gamma^{w,h}(t)&\nonumber\\
    \nonumber\\
    \wb\cdot\nb^w &= \vb\cdot\nb^w \qquad& \text{ on } \gamma^{w,ff}(t)&\nonumber\\
    \wb\cdot\nb^{ff} &= 0 & \text{ on } \gamma^{w,ff}(t)& \label{eq:w-cond-gamma-w-ff}\\
    \wb\cdot(\nb^w\times\nb^{ff}) &= 0 & \text{ on } \gamma^{w,ff}(t)&\nonumber\\
    \nonumber\\
    \wb  &= 0 \qquad& \text{ on }
    \gamma^{b,ff}(t)&.\label{eq:w-cond-gamma-b-ff}
  \end{align}
\end{subequations}
Expressing the boundary conditions~(\ref{eq:w-cond-gamma}) as a given
velocity term $\wb_g$, the evolution equation of $\gamma(t)$ become
\begin{equation}
  \label{eq:evolution-gamma}
  \begin{aligned}
    \frac{\Pa\pale_\gamma}{\Pa t}(\xale,t) &= \wb_g(\pale_\gamma(\xale,t)) \qquad
    &\text{ on } \tilde\gamma&\\
    \pale_\gamma(\xale,0) &= \pale_0(\xale) \qquad &\text{ on } \tilde\gamma&.
  \end{aligned}
\end{equation}

A reconstruction of a reasonable ALE deformation on the entire
$\Gamma(t)$ is then possible by solving an additional elliptic
boundary value problem, coupled with a projection on the surface of
the ship hull and on the free-surface. Given a free-surface
configuration $\eta$ and the deformation $\pale_\gamma$ on the
wireframe $\tilde\gamma$, in order to find a compatible deformation
$\pale$ on the entire $\tilde\Gamma$, we solve the additional problem
\begin{equation}
  \label{eq:laplace-beltrami-on-Gamma}
  \begin{aligned}
    & -\tilde\Delta_s \tilde \gb = -2\tilde\nb \tilde k  \qquad & \text{ on }
    \tilde\Gamma&\\
    & \tilde \gb = \pale_\gamma  & \text{ on } \tilde\gamma&\\
    \\
    &\pale = \mathcal{P}_h\,\tilde \gb  \qquad
    & \text{ on } \tilde\Gamma^h&\\
    & \pale = \mathcal{P}_\eta \,\gale := \gale + (\eta(\gale)-\gale\cdot\eb_z)\eb_z 
    & \text{ on } \tilde\Gamma^w&\\
    \\
    & \pale = \gale 
    & \text{ on } \tilde\Gamma^b\cup\tilde\Gamma^{ff},
  \end{aligned}
\end{equation}
\nomenclature{$\Delta_s$}{Laplace-Beltrami operator}
\nomenclature{$\tilde k(\xale)$}{Mean curvature of $\tilde \Gamma$}
\nomenclature{$\mathcal{P}_h$}{Boat projection operator}
\nomenclature{$\gale$}{Auxiliary domain deformation variable}
where $\tilde k(\xale)$ is the mean curvature of the domain $\tilde\Gamma$,
i.e., the mean curvature of the hull on $\tilde\Gamma^h$ and zero
eveywhere else, while $\mathcal{P}_h$ is a projection operator on the
hull surface. Similarly, $\mathcal P_\eta$ is a (vertical) projection
operator on the free-surface. 

The auxiliary function $\gale$ represents a surface that follows the
waterline deformation $\pale_\gamma$. On the ship hull, it is a
perturbation of the shape of the hull while everywhere else it is a
minimal surface with boundary conditions imposed by $\pale_\gamma$. On
the free-surface, only its $x$ and $y$ components are used to
determine $\pale$, while $\eta$ (which satisfies the kinematic
boundary conditions~(\ref{eq:fsKinematicSemiLag})) imposes the $z$
component.

\subsection{Integro-differential formulation}
\label{sec:integro-diff-form}

Putting everything together, the final integro-differential system is
given by the following problem:

Given initial conditions $\phi_0$ and $\eta_0$ on $\Gamma^w(0)$, and
$\pale_0$ on $\tilde\gamma$, for each time $t\in [0,T]$, find $\pale$,
$\phi$, $\phin$ that satisfy
\begin{subequations}
  \label{eq:continuous-problem}
  \begin{alignat}{2}
    \label{eq:continuous-bie}
    & \int_{\Gamma(t)} \frac{\Pa G}{\Pa n}\phi \d\Gamma-\phi\int_{\Gamma(t)} \frac{\Pa G}{\Pa n} \d\Gamma
    = \int_{\Gamma(t)} \phin G \d\Gamma  & \text{ on } \Gamma(t)\\
    \label{eq:continous-fsDynamicSemiLag}
    & \frac{\delta\phi}{\delta t} = V_\phi(\phi,\phi_n,\eta,\wb)& \qquad  \text{ on } \Gamma^w(t) \\
    &\frac{\delta\eta}{\delta t} = V_\eta(\phi, \phi_n,\eta,\wb)
    & \text{ on }\Gamma^w(t)\label{eq:continuous-fsKinematicSemiLag}\\
    & \phi(\xb, 0) = \phi_0(\xb) &\text{ on }
    \Gamma^w(0)\label{eq:continuous-initial-condition-phi-ALE} \\ 
    &\eta(\xb, 0) = \eta_0(\xb) &\text{ on }
    \Gamma^w(0)\label{eq:continuous-initial-condition-eta-ALE} \\
    \label{eq:continous-non-penetration-condition-hull-ALE}
    &\phin = \overline{\phin} & \text{ on } \Gamma^{N}(t) \\
    \label{eq:continuous-ALE-velocity-condition}
    &\wb\cdot\nb = 0 & \text{ on } \Gamma^{N}(t) \\
    &-\tilde\Delta_s\gale = -2\tilde\nb k \qquad &\text{ on }
    \tilde\Gamma\setminus\tilde\gamma\\
    &\frac{\Pa\gale_\gamma}{\Pa t}(\xale,t) = \wb_g(\gale(\xale,t)) \qquad
    &\text{ on } \tilde\gamma\\
    &\pale(\gale,0) = \pale_0(\gale) \qquad &\text{ on }
    \tilde\gamma&\\
    & \pale = \mathcal P\,\gale &\text{ on }
    \tilde\Gamma.
  \end{alignat}
\end{subequations}
where we used the shorthand notations
\begin{subequations}
  \begin{align}
    \wb &:= \frac{\delta \pb}{\delta t}\\
    \label{eq:def-total-V}
    \vb & := \Vb_\infty-\Vb_\sys+\nablab\phi\\
    \label{eq:def-V-phi}
    V_\phi(\phi,\phi_n,\eta,\wb) &:= (\wb-\vb)\cdot\nablab\phi-g\eta + \ab_{\infty}\cdot\xb+
    \frac{1}{2}|\nablab\phi|^2 -\mu\phin \\
    \label{eq:def-V-eta}
    V_\eta(\phi,\phi_n,\eta,\wb) &:= (\wb-\vb)\cdot\nablab\eta + \vb\cdot
    \eb_z \\
    \Gamma^N(t) &:= \Gamma^{h}(t)\cup\Gamma^{b}(t)\cup\Gamma^{ff}(t) \\
    \overline\phin &:=
    \begin{cases}
      (\Vb_\boat-\Vb_\sys)\cdot\nb & \text{ on } \Gamma^h(t) \\
      0 & \text{ on } \Gamma^{b}(t)\cup\Gamma^{ff}(t),
    \end{cases}\\
    \label{eq:projection-operator}
    \mathcal P \gale& :=
    \begin{cases}
      \mathcal{P}_h \,\gale  & \text{ on } \tilde\Gamma^h\\
      \mathcal{P}_\eta\, \gale:= \gale + (\eta(\gale,t)-\gale\cdot\eb_z)\eb_z 
       & \text{ on } \tilde\Gamma^w\\
      \gale  & \text{ on } \tilde\Gamma^b\cup \tilde\Gamma^{ff},
    \end{cases}
  \end{align}
\end{subequations}
and both the potential and the pressure in the entire domain can be
obtained by postprocessing the solution to
Problem~(\ref{eq:continuous-problem}) with the boundary integral
representation~(\ref{eq:bie}) and with Bernoulli's
Equation~(\ref{eq:bernoulli}). 

\nomenclature{$V_\phi$}{Right hand side for ALE derivative of the
  perturbation potential $\phi$}
\nomenclature{$V_\eta$}{Right hand side for ALE derivative of the
  surface elevation $\eta$}
\nomenclature{$\mathcal{P}_\eta$}{Projection operator on the free-surface}

The full gradient of the perturbation potential on the surface
$\Gamma(t)$ that appears in Equations~(\ref{eq:def-total-V}),
(\ref{eq:def-V-phi}) and~(\ref{eq:def-V-eta}) is constructed from the
surface gradient of $\phi$ and from the normal gradient $\phin$ as
\begin{equation}
  \nablab \phi(\xb,t) := \nablab_s\phi(\xb,t) +\phin(\xb,t)\nb.
\end{equation}

A numerical discretization of the continuous
Problem~(\ref{eq:continuous-problem}) is done on the \emph{fixed}
boundary $\tilde\Gamma$ of the reference domain $\Oale$, with
independent variable $\xale$ which will label node locations in a
reference computational grid, and the motion $\pale(\xale,t)$ will
denote the trajectory of the computational nodes.


\section{Numerical discretization}
\label{sec:numer-discr}

To approximate the continuous problem, we introduce a decomposition
$\tilde\Gamma_{h}$ of the domain boundary $\tilde\Gamma$. Such
partition is composed of three-dimensional quadrilateral cells
(indicated with the index $K$), which satisfy the following regularity
assumptions:
\begin{enumerate}
\item $\overline{\tilde\Gamma} = \cup \{ K \in \tilde\Gamma_{h} \}$;

\item
Any two cells $K,K'$ only intersect in common faces, edges, or vertices;

\item
The decompositions $\tilde\Gamma_{h}$ matches the decomposition
$\tilde\Gamma = \tilde\Gamma^w \cup \tilde\Gamma^h \cup\tilde\Gamma^b
\cup \tilde\Gamma^{ff}$.
\end{enumerate}
On the decomposition $\tilde\Gamma_{h}$, we look for solutions
$\displaystyle(\pb_h,\phi_h, \phin_h)$ in the finite dimensional
spaces ${Y}_{h}$, ${V}_{h}$, and ${Q}_{h}$ defined as
\begin{alignat}{5}
\label{eq:space-Y-h}
{Y}_{h} &:= \Bigl\{ \ub_h \in \mathcal{C}^0(\tilde\Gamma_h)^3 \,&&\big|\, \ub_{h|K} &&\in \mathcal{P}_Y(K)^3, \, K &&\in \tilde\Gamma_h \Bigr\} &&\equiv \vssp\{ \vb_{h}^{i} \}_{i=1}^{N_{Y}}\\
\label{eq:space-phi-h}
{V}_h &:= \Bigl\{ \phi_h \in \mathcal{C}^0(\tilde\Gamma_h)  \,&&\big|\, \phi_{h|K}  &&\in \mathcal{P}_V(K), \, K &&\in \tilde\Gamma_h \Bigr\} &&\equiv \vssp\{ \varphi_{h}^{i} \}_{i=1}^{N_{V}}
\\
\label{eq:space-phin-h}
{Q}_h &:= \Bigl\{ \gamma_h \in \mathcal{C}^0(\tilde\Gamma_h)  \,&&\big|\, \gamma_{h|K}  &&\in \mathcal{P}_Q(K), \, K &&\in \tilde\Gamma_h \Bigr\} &&\equiv \vssp\{ \tau_{h}^{i} \}_{i=1}^{N_{Q}}.
\end{alignat}
Here, $\mathcal{C}^0(\tilde\Gamma_h)^d$ is the space of continuous
functions of $d$ components over $\tilde\Gamma_h$. 
$\mathcal{P}_{Y}(K)$, $\mathcal{P}_{V}(K)$ and $\mathcal{P}_{Q}(K)$
indicate the polynomial spaces of degree $r_{Y}$, $r_{V}$ and $r_{Q}$
respectively, defined on the cells $K$. Finally, $N_Y$, $N_V$ and
$N_Q$ denote the dimensions of each finite dimensional space.

\nomenclature{$K$}{A quadrilateral cell of the discretization}
\nomenclature{$\tilde\Gamma_h$}{The discretization of the surface
  $\Gamma$}
\nomenclature{$\mathcal{P}_{Y}(K)$}{Polynomial space of degree $r_{Y}$}
\nomenclature{$\mathcal{P}_{Q}(K)$}{Polynomial space of degree $r_{Q}$}
\nomenclature{$\mathcal{P}_{V}(K)$}{Polynomial space of degree
  $r_{V}$}
\nomenclature{$Y_{h}$}{Finite dimensional space of deformations}
\nomenclature{$V_h$}{Finite dimensional space of potential}
\nomenclature{$Q_h$}{Finite dimensional space of normal derivative of
  the potential}
\nomenclature{$\vb_{h}$}{Basis functions for $Y_h$}
\nomenclature{$\varphi_{h}$}{Basis functions for $V_h$}
\nomenclature{$\tau_{h}$}{Basis functions for $Q_h$}
\nomenclature{$N_V$}{Dimension of the space $V_h$}
\nomenclature{$N_Y$}{Dimension of the space $Y_h$}
\nomenclature{$N_Q$}{Dimension of the space $Q_h$}

The most common approach for the solution of the boundary integral
equation~(\ref{eq:bie}) in the engineering community is the
\emph{collocation boundary element method}, in which the continuous
functions $\phi$ and $\nablab \phi \cdot \nb$ are replaced by their
discrete counterparts and the boundary integral equation is imposed
at a suitable number of points on $\Gamma(t)$.

Once a geometric representation of the reference domain
$\tilde\Gamma_h$ is available, we could in principle employ arbitrary
and independent discretizations for each of the functional
spaces $V_h$, $Q_h$ and $Y_h$. We choose instead to adopt an
\emph{iso-parametric} representation, in which
the same parametrization is used to describe both the surface geometry and
the physical variables. Thus, the deformed surface $\Gamma_h(t)$
(i.e., the map $\pale(\xale,t)$),
the surface potential $\phi_h(\xb,t)$ and its normal derivative
$\displaystyle\phin(\xb,t)$, are discretized on the $k$-th panel as
\begin{equation}
  \label{eq:space-choice}
  Q_h = V_h, \qquad Y_h = V_h^3 = \vssp\{ \varphi_{h}^{i}\eb_x, 
  \varphi_{h}^{i}\eb_y,\varphi_{h}^{i}\eb_z \}_{i=1}^{N_{V}},
\end{equation}
where $\eb_{x,y,z}$ are unit basis vectors identifying the directions
of the $x$, $y$ and $z$ axis.
We indicate with the notation $\{\phi\}$ and
$\displaystyle\left\{\phin\right\}$ the column
vectors of time-dependent coefficients $\phi^i(t)$ and
$\displaystyle\phin^j(t)$ such that
\begin{equation}
  \label{eq:coefficient-notation}
  \begin{split}
    \phi_h(\xb, t) & := \sum_{i=1}^{N_V}
    \phi^i(t)\varphi^i_h(\pale_h^{-1}(\xb,t)) \qquad\text{ on }\Gamma_h(t)\\
    \phin_h(\xb, t) &:= \sum_{j=1}^{N_V}
    \phin^j(t)\varphi^j_h(\pale_h^{-1}(\xb,t))  \qquad\text{ on }\Gamma_h(t).
  \end{split}
\end{equation}
The map $\pale_h^{-1}(\xb,t)$ is the inverse of the discretizaed ALE
deformation, which reads
\begin{equation}
  \label{eq:control-points}
  \pale_h(\xale, t) := \sum_{k=0}^{N_V}\xb^k(t)\varphi^k_h(\xale) \qquad\text{ on }\tilde\Gamma,
\end{equation}
where $\xb^k$ represents the time-dependent location of the vertices or
control points that define the current configuration of
$\Gamma(t)$. 

To distinguish matrices from column vectors, we will indicate matrices
with the bracket notation, e.g.., $[M]$. The ALE derivatives of the
finite dimensional $\phi_h(\xb,t)$ and $\pb_h(\xb,t)$ are
time parametrized finite dimensional vectors in
$V_h$ and $Y_h$, reading
\begin{equation}
  \label{eq:ale-derivative-discrete}
  \begin{split}
    \frac{\delta \phi_h}{\delta t}(\xb, t) & = \sum_{i=1}^{N_V}
    \frac{\Pa
      \phi^i}{\Pa t}(t)\varphi^i_h(\pale_h^{-1}(\xb,t)) \\
    \wb_h(\xb, t) := \frac{\delta \pb_h}{\delta t}(\xb, t) & = \sum_{i=1}^{N_V} \frac{\Pa
      \xb^i}{\Pa t}(t)\varphi^i_h(\pale_h^{-1}(\xb,t)).
  \end{split}
\end{equation}
Here, each function is identified by the coefficients $\{\phi\}'$ and by the
control points $\{\xb\}'$, where the $'$ denotes derivation in time.

\nomenclature{$[M]$}{Mass matrix}

Finally, we can reconstruct the full discrete gradient $\nablab\phi_h$ on
$\Gamma_h(t)$ using the discrete version of the surface gradient
$\nablab_s$ and the normal gradient $\phin_h$:
\begin{equation}
  \label{eq:surface-discrete-gradient}
  \begin{split}
    \nablab_s\phi_h(\xb, t) & := \sum_{i=1}^{N_V} \phi^i(t)
    \nablab_s\varphi^i_h(\pale_h^{-1}(\xb,t)) \qquad\text{ on
    }\Gamma_h(t)\\
    \nablab\phi_h(\xb, t) & = \nablab_s\phi_h(\xb, t) + \phin_h(\xb,t)\nb  \qquad\text{ on }\Gamma_h(t).
  \end{split}
\end{equation}

\subsection{Collocation boundary element method}
\label{sec:collocation-bem}

The discrete version of the BIE, written for an arbitrary point $\yb$ on the domain boundary,
reads

\begin{multline}
    \alpha(\yb,t)\phi(\yb,t)  = \\
    -\sum_{k=1}^{M}\sum_{i=1}^{N_V} \phi^i(t)
    \int_{\hat K} \frac{\Pa G}{\Pa n}(\yb-\xb^k(u,v,t))
    \varphi^i_k(u,v)  J^k(u,v,t) \d u \d v \\
    + \sum_{k=1}^{M}\sum_{i=1}^{N_V} \left(\frac{\Pa \phi}{\Pa
        n}(t)\right)^i \int_{\hat K} G(\yb-\xb^k(u,v,t))  \varphi^i_k(u,v)  J^k(u,v,t) \d
    u \d v.
\label{eq:discBIE1}
\end{multline}
Here, we made use of the iso-parametric representation described in~\ref{sec:iso-param},
to decompose the integrals into the local contributions of the $N_L$ basis functions
in each of the $M$ panels of the triangulation.

The numerical evaluation of the panel integrals appearing in equation
(\ref{eq:discBIE1}) needs some special treatment, due to the presence
of the singular kernels $G(\yb-\xb)$ and $\frac{\Pa G}{\Pa
  n}(\yb-\xb)$. Whenever $\yb$ is not a node of the integration panel,
the integral argument is not singular, and standard Gauss quadrature
formulas are used. If $\yb$ is a node of the integration panel, the
integral kernel is singular and special quadrature rules are used,
which remove the singularity by performing an additional change of
variables (see, for example, \cite{lachatWatson}). In the framework of
collocated BEM, an alternative possibility would have been represented by
desingularized methods, in which the fundamental solutions are centered
at points that are different from the evaluation points. Typically, this
is obtained by centering the fundamental solutions at points that are
slightly outside the domain. Although these methods avoid dealing with
singular integrals, they pose problems on establishing a general rule
for suitable positioning of the fundamental solutions centers.
In the case at hand, the domain presents sharp edges and narrow corners
(typically found at the bow or stern of a hull) which make the latter
task nontrivial.

Writing a boundary integral equation for each node $\xb^i$, $i=1,\dots,N_V$ of
the computational grid, we finally obtain the system
	
\begin{equation}
\label{eq:solvingSysFMM}
\left[\alpha\right]\left\{\phi\right\} + \left[N\right]\left\{\phi\right\} =
\left[D\right]\left\{\phin\right\},
\end{equation}
where
	
\begin{itemize}
	
\item $\left\{\phi\right\}$ and $\left\{\phin\right\}$
  are the vectors containing the potential and its normal
  derivative node values, respectively;
  
\item $\left[\alpha\right]$ is a diagonal matrix composed by the
  $\alpha(\xb_i(t))$ coefficients;
  
\item $\left[D\right]$ and $\left[N\right]$ are the Dirichlet and
  Neumann matrices respectively whose elements are
\begin{eqnarray}
D_{ij} &=& \sum_{k=1}^{M}
                       \int_{\hat K}
		       G(\xb_i(t)-\xb^k(u,v,t))
                       \varphi^j_k(u,v)J^k(u,v,t) \d u \d v \\
N_{ij} &=& \sum_{k=1}^{M}
                       \int_{\hat K}
		       \frac{\Pa G}{\Pa n}(\xb_i(t)-\xb(u,v,t))
		         \varphi^j_k(u,v)J^k(u,v,t) \d u \d v. \\
\end{eqnarray}
\end{itemize}	

\nomenclature{$\left\{\phi\right\}$}{Vector of potential node values}
\nomenclature{$\left\{\phin\right\}$}{Vector of normal derivative node values}
\nomenclature{$\left[\alpha\right]$}{Diagonal matrix composed by the
  $\alpha(\xb_i(t))$ coefficients}
\nomenclature{$\left[D\right]$}{Dirichlet matrix}
\nomenclature{$\left[N\right]$}{Neumann matrix}

The evaluation of the nodal values for the solid angle fractions
$\alpha_i$ appearing in the BIE equation is obtained considering the
solution to Laplace problem~(\ref{eq:laplace-problem}) when
$\phi\equiv 1$ in $\Omega(t)$. In such case, system
(\ref{eq:solvingSysFMM}) reads

\begin{equation}  
\left[\alpha\right]\left\{1\right\} +
\left[N\right]\left\{1\right\} = 0,
\end{equation}
which implies
\begin{equation} \label{eq:rigidModes} 
\alpha_i = -\sum_{j=1}^{N} N_{ij} \qquad i = 1,\dots,N.
\end{equation}

This technique is usually referred to as the Rigid Mode Technique, or
Rigid Body Mode (RBM) (see, for example,\cite{Cruse1974}). It can be
interpreted as an indirect regularization of the Neumann matrix $[N]$,
which ensures that the matrix $[\alpha]+[N]$ possesses a zeroth
eigenvalue corresponding to the rigid mode of the system.

\begin{figure}[htb!]        
\centerline{
  \ifpdf
  \resizebox{.5\textwidth}{!}{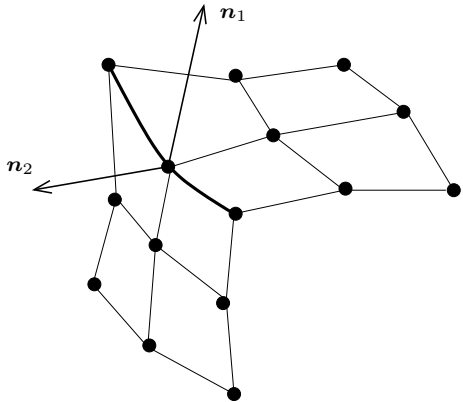}
  \else
  \resizebox{.5\textwidth}{!}{\input{figures/edgesTreat.pstex_t}}
  \fi
}
\caption{The two different normal unit vectors of a node placed on  
         an edge \label{fig:edge}}	
\end{figure}

It is worth noting here, that there are arguments in favor of selecting
discontinuous discrete representations for $\phin = \nablab\phi\cdot\nb$.
Such considerations are based on the fact that surface normals are not
necessarily continuous across neighboring panels. In most cases though,
the nodal $\{\phin\}$ values obtained through the solution of
system (\ref{eq:solvingSysFMM}) represent a very reasonable approximation
of the continuous normal potential gradient. But whenever the domain boundary
presents sharp features, such as corners or edges (see Fig. \ref{fig:edge}),
a continuous approximation of $\{\phin\}$ given by
system (\ref{eq:solvingSysFMM}) becomes extremely inaccurate, the exact
normal potential gradient itself being not continuous.

To overcome such problem, in this work we employ a technique for the treatment
of edges and corners, which was first developed by Grilli and
Svendsen (\cite{grilliCorner90}). In this framework, the computational grid
is first prepared so that on each edge separating two different boundary condition
zones, the mesh nodes are duplicated. Hence, two distinct nodes
are present in the mesh, each belonging to one of the two elements on the edge,
and each having a different normal unit vector. In this way, the number of degrees
of freedom of the system has been increased so as to allow discontinuous edge
values for the normal derivative of the potential. 

\subsection{SUPG stabilization}
\label{sec:semi-discrete}

For the time evolution of the dynamic and kinematic boundary
conditions~(\ref{eq:continous-fsDynamicSemiLag})
and~(\ref{eq:continuous-fsKinematicSemiLag}), we need to evaluate the
surface gradient of the veocity potential, as it appears in
Equation~(\ref{eq:surface-discrete-gradient}).

The gradient of the perturbation potential is not, in general,
continuous across the edges of the panels that compose $\Gamma_h(t)$.
As a consequence, the right hand side of
equations~(\ref{eq:continous-fsDynamicSemiLag}) is not single valued
at the location of the computational grid nodes. Thus, it is not possible
to write a direct evolution equation for the nodes of the computational
grid and for the  corresponding potential values.

A possible solution to this problem would be collocating
the time evolution equations at marker
points placed on the internal surface of each cell (\cite{KjelContJans-2011}).
Although this strategy would result in a single valued right hand side of
equations~(\ref{eq:continous-fsDynamicSemiLag}), it would require
that a new computational mesh is reconstructed from the updated
markers position at each time step. A different approach
consists in using using smooth finite dimensional spaces,
as in~\cite{grilli2001}. This choice would allow
for the collocation of the evolution equations at the mesh nodes,
but would in turn require the use of high order panels.

In this work, we choose instead to impose the evolutionary
boundary conditions in weak form. Specifically, we employ an $L^2$
projection in the $V_h$ space, to evaluate the right hand side of
equations~(\ref{eq:continous-fsDynamicSemiLag})
and~(\ref{eq:continuous-fsKinematicSemiLag}), namely
\begin{subequations}
  \label{eq:l2projection}
  \begin{alignat}{2}
     & \left(\frac{\delta\phi}{\delta t},
       \varphi\right)_w = \left(V_\phi, \varphi\right)_w\qquad
     &\forall \varphi \in V_h \label{eq:continuous-fsDinamicSemiLag-weak}\\
     &\left(\frac{\delta\eta}{\delta t},\varphi\right)_w = \big(V_\eta,\varphi\big)_w
     &  \forall \varphi \in V_h \label{eq:continuous-fsKinematicSemiLag-weak},
  \end{alignat}
\end{subequations}
where 
\begin{equation}
  \begin{split}
    \left(a,b\right)_w &= \int_{\Gamma^w(t)}ab\ \d\Gamma\\
    \left(a,b\right) &= \int_{\Gamma(t)}ab\ \d\Gamma.
  \end{split}
\end{equation}

\begin{figure}[htb!]
\begin{center}
\includegraphics[width=\textwidth]{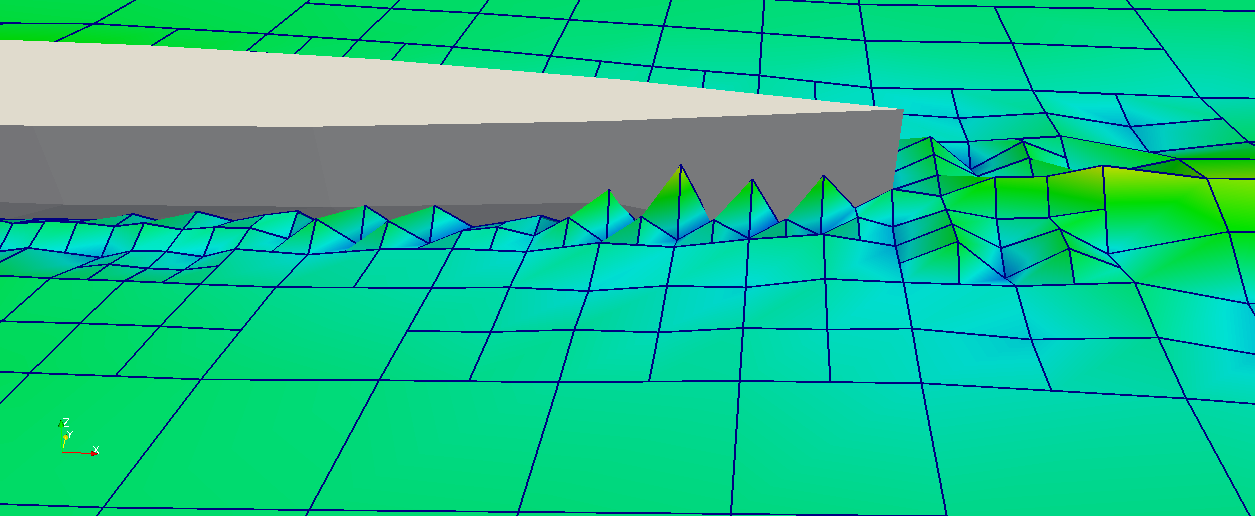}
\end{center}
\caption{An example of the sawtooth instability developing on the
         stern of the hull without stabilization.
         \label{fig:sawtooth}}
\end{figure}

This approach leads to extremely accurate estimations of the forcing terms
in the free-surface evolution equations~(\ref{eq:def-V-phi})
and~(\ref{eq:def-V-eta}), with a very small computational cost.
Unfortunately, the equations right hand sides
both contain transport terms (respectively $\nablab\eta\cdot\left(\wb
  -\nablab\phi-\Vb_\infty+\Vb_\sys\right)$ and
$\nablab\phi\cdot\left(\wb -\nablab\phi-\Vb_\infty+\Vb_\sys\right)$),
that become dominant whenever $(\Vb_\sys-\Vb_\infty)$ is very
different from $\nablab\phi$. This causes a sawtooth numerical instability
which in most cases develops in proximity of the hull stern, with
consequent blow up of the simulations (an example of such instability
is given in Figure~\ref{fig:sawtooth}).


A natural, inexpensive and consistent stabilization algorithm which is able toreduce
the observed instabilities is the
Streamwise Upwind Petrov--Galerkin (SUPG) scheme (see, for
example,~\cite{hughes1979,tezduyar2003computation}). The SUPG stabilization
consists in replacing the plain $L^2$ projection in
system~(\ref{eq:l2projection}) with the weighted projection
\begin{subequations}
  \label{eq:supgl2projection}
  \begin{alignat}{2}
    & \left(\frac{\delta\phi}{\delta t},
      \varphi+\db\cdot\nablab_s\varphi\right)_w = \left(V_\phi,
      \varphi+\db\cdot\nablab_s\varphi\right)_w\qquad
    &\forall \varphi \in V_h \label{eq:continuous-fsDinamicSemiLag-weak-supg}\\
    &\left(\frac{\delta\eta}{\delta
        t},\varphi+\db\cdot\nablab_s\varphi\right)_w =
    \left(V_\eta,\varphi+\db\cdot\nablab_s\varphi\right)_w & \forall
    \varphi \in V_h\label{eq:continuous-fsKinematicSemiLag-weak-supg},
  \end{alignat}
\end{subequations}
where 
\begin{equation}
  \label{eq:definition-db}
  \db := \tau \left(\frac{\vb-\wb}{|\vb-\wb|}\right).
\end{equation}
$\tau$ is a positive stabilization parameter which involves a
measure of the local length scale (i.e. the ``element length'') and
the local Reynolds and Courant numbers. Element lengths and
stabilization parameters were proposed for the SUPG formulation of
incompressible and compressible flows in~\cite{hughes1979}, and an in
depth study of the stabilization properties for free boundary problems
was presented in~\cite{tezduyar2003computation}. However, to the best
of the authors' knowledge, this is the first time that such
stabilization technique is applied directly to a boundary value
problem defined on a curved surface.

Expressing system~\ref{eq:supgl2projection} in matrix form, we get
\begin{subequations}
  \begin{align}
    &[M][P_{\Gamma^w}]\{\phi\}' = [P_{\Gamma^w}]\{V_\phi\} \\
    &[M][P_{\Gamma^w}]\{\eta\}' = [P_{\Gamma^w}]\{V_\eta\},
  \end{align}
\end{subequations}
where the vectors and ---sparse--- matrices elements are computed by
\begin{subequations}
\begin{align}
M^{ij} := &
\left(\varphi^j,\varphi^i+\db\cdot\nablab_s\varphi^i\right) \nonumber\\
 = &\sum_{k=1}^{M}
\int_{\hat K}
\varphi^j_k(u,v)\big(\varphi^i(u,v)+\db\cdot\nablab_s\varphi^i(u,v)\big)J^k(u,v,t)
\d u \d v\\
P_{A}^{ij} := &
\begin{cases}
  \delta_{ij} & \text{ if } \xb^i(t) \in A(t) \\
  0 & \text{ otherwise}
\end{cases}\\
V_\phi^i := & \left(V_\phi,\varphi^i+\db\cdot\nablab_s\varphi^i\right)
\nonumber\\
= & \sum_{k=1}^{M}
\int_{\hat K}
V_\phi(u,v)\big(\varphi^i(u,v)+\db\cdot\nablab_s\varphi^i(u,v)\big)J^k(u,v,t)
\d u \d v\\
V_\eta^i := & \left(V_\eta,\varphi^i+\db\cdot\nablab_s\varphi^i\right)
\nonumber\\
= & \sum_{k=1}^{M}
\int_{\hat K}
V_\eta(u,v)\big(\varphi^i(u,v)+\db\cdot\nablab_s\varphi^i(u,v)\big)J^k(u,v,t) \d u \d v.
\end{align}
\end{subequations}

\subsection{Semi-discrete smoothing operator}

The semi-discrete version of the smoothing
problem~(\ref{eq:laplace-beltrami-on-Gamma}) can be obtained with a
finite element implementation of the scalar Laplace-Beltrami operator,
and its application to the different components of $\pale$.

The weak form of the Laplace-Beltrami operator on $\tilde\Gamma$ of a
scalar function $u$ in $V$ with Dirichlet boundary conditions $u_g$ on
$\tilde\gamma$ and forcing term $f$ reads
\begin{equation}
  \label{eq:weak-laplace-beltrami} 
  \begin{aligned}
    &\left(\tilde\nablab_s u,
      \tilde\nablab_s\varphi\right)_{\tilde\Gamma} = \left(f,
      \varphi\right)_{\tilde\Gamma}  \qquad &\forall
    \varphi \in V_0\\
    & u = u_g \qquad& \text{ on } \tilde\gamma.
  \end{aligned}
\end{equation}
Here, $V_0$ denotes the space of functions $\varphi$ in $V$
such that their trace on $\tilde\gamma$ is zero (see, for example,
\cite{BonitoNochettoPauletti-2010-a} and the references therein for
some details on the numerical implementation of the Laplace-Beltrami
operator). 

The semi-discrete form of Equation~(\ref{eq:weak-laplace-beltrami})
can be written as
\begin{equation}
  [K]\{u\} = \{F\}, \label{eq:semidiscrete-laplace-beltrami}
\end{equation}
where
\begin{equation}
  \begin{aligned}
    K^{ij} := &
    \left(\tilde\nablab_s\varphi_h^j,\tilde\nablab_s\varphi_h^i\right) \\
    = &\sum_{k=1}^{M} \int_{\hat K} \tilde\nablab_s\varphi_k^j(u,v)
    \cdot \tilde
    \nablab_s \varphi_k^i(u,v) \tilde J_k(u,v) \d u \d v\\
    F^{i} := &
    \left(f,\varphi_h^i\right) \\
    = &\sum_{k=1}^{M} \int_{\hat K}
    f_k^j(u,v) \varphi_k^i(u,v) \tilde J_k(u,v) \d u \d v.
  \end{aligned}
\end{equation}

The discrete Laplace-Beltrami operator is solved for the auxiliary
vector variable $\gale$, whose finite dimensional representation is
given by $\{\gb\}$. The forcing terms in the system are given by
the mean curvature along the normal of $\tilde\Gamma$. 
In this case we write
\begin{equation}
  [\Kb] \{\gb\} = \{\kb\},
\end{equation}
where
\begin{equation}
  \begin{split} 
    [\Kb] &:= 
     \begin{bmatrix}
       [K] & 0 & 0 \\
       0 & [K] & 0\\
       0 & 0 & [K]
     \end{bmatrix}\\
     \{\kb\} &:= - 2 k
     \begin{Bmatrix}
       \nb_x \\
       \nb_y \\
       \nb_z
     \end{Bmatrix}.
  \end{split}
\end{equation}
and the full matrix form of
Problem~(\ref{eq:continuous-problem}) finally reads
\begin{subequations}
  \label{eq:semi-discrete-problem}
  \begin{align}
    \label{eq:semidiscrete-bem}
    &\left[\alpha\right]\left\{\phi\right\} + \left[N\right]\left\{\phi\right\} -
    \left[D\right]\left\{\phin\right\}  &= 0 \\
    &[M][P_{\Gamma^w}]\{\phi\}' - [P_{\Gamma^w}]\{V_\phi\}& = 0 \label{eq:semidiscrete-bem-eta}\\
    &[M][P_{\Gamma^w}]\{\eta\}' - [P_{\Gamma^w}]\{V_\eta\} &= 0 \label{eq:semidiscrete-bem-phi}\\
    & [P_{\Gamma^w}]\{\phi (0)\} - [P_{\Gamma^w}]\{\phi_0\}& = 0\\ 
    &[P_{\Gamma^w}]\{\eta (0)\} - [P_{\Gamma^w}]\{\eta_0\} &= 0\\
    &[I-P_{\Gamma^w}]\{\phin \} - [I-P_{\Gamma^w}]\overline{\{\phin\}}
    &= 0\\
    &[P_{\gamma}]\{\gb\}' -[P_{\gamma}]\{\wb_g\} & =0\\
    &[\Kb]\{\gb\} - \{\kb\} &= 0 \\
    &[\mathcal P]\{\gb\} - \{\xb\} &= 0,
  \end{align}
\end{subequations}
where $\mathcal P$ is a numerical implementation of the projection
operator~(\ref{eq:projection-operator}). On the hull surface, this
operator is easily implemented analytically for simple
hull shapes, such as that of the Wigley hull. In a more general case,
it is desirable to have an implementation of the projection operator
which can directly interrogate the CAD files describing the hull
surface. In this work, we present results obtained with the former
projection. A CAD based projection, which makes use of 
the OpenCASCADE library~\cite{opencascade}, is currently being
implemented. Some work is still required though, to
render our full discretization robust with respect to arbitrary hull
geometries.





\subsection{Time discretization}
\label{sec:time-discretization}

System~(\ref{eq:semi-discrete-problem}) can be recast in the following
form

\begin{equation}
\label{eq:dae_formulation}
  F(t, y, y') = 0,
\end{equation}
where we grouped the variables of the system in the vector $y$:
\begin{equation}
\label{eq:dae_unknowns}
y = \left\{\begin{array}{c}
           \left\{\xb\right\} \\
           \left\{ \phi\right\} \\
           \left\{ \phin\right\}
           \end{array}\right\}.
\end{equation}

Equation (\ref{eq:dae_formulation}) represents a system of nonlinear differential algebraic
equations ({DAE}), which we solve using the {IDA} package of the {SUNDIALS} OpenSource
library~\cite{sundials2005}.  As stated in the package
documentation (see p.~374 and~375 in \cite{sundials2005}):\footnote{We quoted
directly from the {SUNDIALS} documentation. However, we adjusted the notation so as to
be consistent with ours and we numbered equations according to their order in this paper.}
\begin{quote}
The integration method in {IDA} is variable-order,
variable-coefficient {BDF} [backward difference formula], in fixed-leading-coefficient form.
The method order ranges from 1 to 5, with the {BDF} of order $q$ given by the multistep formula
  \begin{equation}
    \label{eq:bdf_of_order q}
    \sum_{i=0}^{q} \alpha_{n,i} y_{n-i} = h_n \dot{y}_n,
  \end{equation}
where $y_{n}$ and $\dot{y}_{n}$ are the computed approximations to $y(t_{n})$ and $y'(t_{n})$,
respectively, and the step size is $h_{n} = t_{n} - t_{n-1}$. The coefficients $\alpha_{n,i}$ are
uniquely determined by the order $q$, and the history of the step sizes.  The application of the
{BDF} [in Eq.~\eqref{eq:bdf_of_order q}] to the {DAE} system [in Eq.~\eqref{eq:dae_formulation}]
results in a nonlinear algebraic system to be solved at each step:
  \begin{equation}
    \label{eq:dae_algebraic_system}
    R(y_{n}) \equiv F\biggl(t_{n}, y_{n},  h_{n}^{-1} \sum_{i=0}^{q}
    \alpha_{n,i} y_{n-i}\biggr) = 0.
  \end{equation}
Regardless of the method options, the solution of the nonlinear
system [in Eq.~\eqref{eq:dae_algebraic_system}] is accomplished with some form of Newton
iteration. This leads to a linear system for each Newton correction, of the form
  \begin{equation}
    \label{eq:dae newton correction}
    J[y_{n,m+1}-y_{n,m}] = -R(y_{n,m}),
  \end{equation}
where $y_{n,m}$ is the $m$th approximation to $y_{m}$. Here $J$ is  some approximation to the system Jacobian
\begin{equation}
  \label{eq:dae_Jacobian}
  J = \frac{\partial R}{\partial y} = \frac{\partial F}{\partial y} +\alpha \frac{\partial F}{\partial y'},
\end{equation}
where $\alpha = \alpha_{n,0}/h_{n}$. The scalar $\alpha$ changes whenever the step size or method order changes. 
\end{quote}
In our implementation, we assemble the residual $R(y_{n,m})$ at each
Newton correction, and let SUNDIALS compute an approximation of the Jacobian
in Eq.~(\ref{eq:dae_Jacobian}). The final system is solved using a preconditioned GMRES iterative
method (see, e.g., \cite{GolubVan-Loan-1996-a}). 

Despite the increase in computational cost due the implicit solution scheme, the DAE system
approach presents several advantages with respect to explicit resolution techniques. 
First, it is worth pointing out that among the unknowns in equation~(\ref{eq:dae_unknowns}), the
coordinates of all the grid nodes (except for the vertical coordinates of free-surface nodes)
appear in the DAE system as \emph{algebraic} components. Their evolution is in fact not described by a
differential equation, but computed through the smoothing operator. Yet, the time derivative of
such coordinates, \emph{i.e.}: the ALE velocity $\wb$ is readily available through the evaluation
of the BDF (\ref{eq:bdf_of_order q}). Such velocites are used in the differential
Equations (\ref{eq:semidiscrete-bem-eta}) and (\ref{eq:semidiscrete-bem-phi}) which appear in the
DAE system. In particular, at each time step, the convergence of Newton
corrections (\ref{eq:dae newton correction}) ensures that the vertical velocity of the nodes is
the one corresponding to the $\wb$ velocity originated by the horizontal nodes displacement computed
by the smoothing operator. In a similar fashion, the DAE solution algorithm computes the ALE time
derivative of the velocity potential at each Newton correction. Such derivative is plugged into Bernoulli's
equation~(\ref{eq:bernoulli}) to evaluate the pressure on the whole domain boundary $\Gamma(t)$,
without requiring the solution of additional boundary value problems for $\frac{\Pa \phi}{\Pa t}$.
Finally, the resulting pressure field is integrated on the ship wet surface $\Gamma^h(t)$ to obtain the
pressure force acting on the ship. Since this operation is done at the level of each Newton
correction, possible rigid motions of a hull along its six degrees of freedom can be accounted for
in a very natural way in the DAE framework, by adding the six differential equations of
motion governing the unknown rigid displacements to the DAE system. The latter ---strongly coupled---
fluid structure interaction model is currently under development and results will be presented in
future contributions.

\subsection{Adaptive mesh refinement}

There are two main advantages of using a Galerkin formulation for the
evolution equation of the free-surface, as well as for the computation of the
full ALE deformation on the surface $\Gamma(t)$. On one
hand, it allows the use of fully unstructured meshes, with an immediate
simplification in the mesh generation. On the other hand, it 
enables the use of a wide set of local refinement strategies based
on a posteriori error estimators, which are rather popular in the
finite element community. As a result, the mesh generation and adaptation
are fully automated, and the resulting computational grids
are shaped based on the characteristic of the
problem solution, rather than on a-priori heuristic choices.

In this work, we use a modification of the gradient recovery error
estimator by Kelly, Gago, Zienkiewicz and
Babuska~\cite{kelly1983posteriori,gago1983posteriori}, a choice mostly
motivated by its simplicity (see~\cite{AinsworthOden-1997-a} for more
details on this and other error estimators).

\begin{figure}[htb!]        
\begin{tabular}{c c}
\includegraphics[width=.48\linewidth]{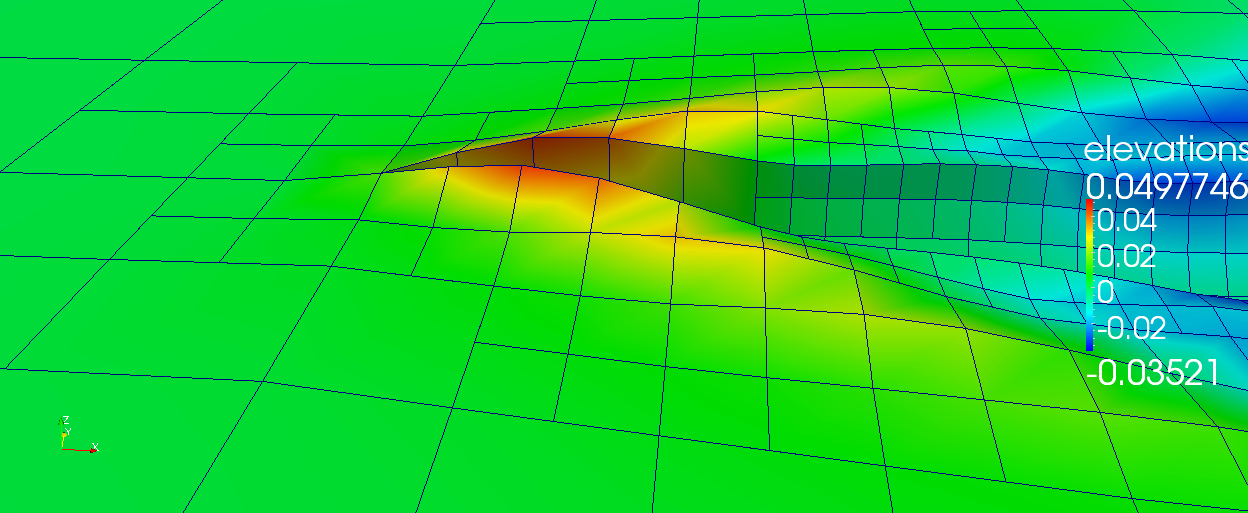} &
\includegraphics[width=.48\linewidth]{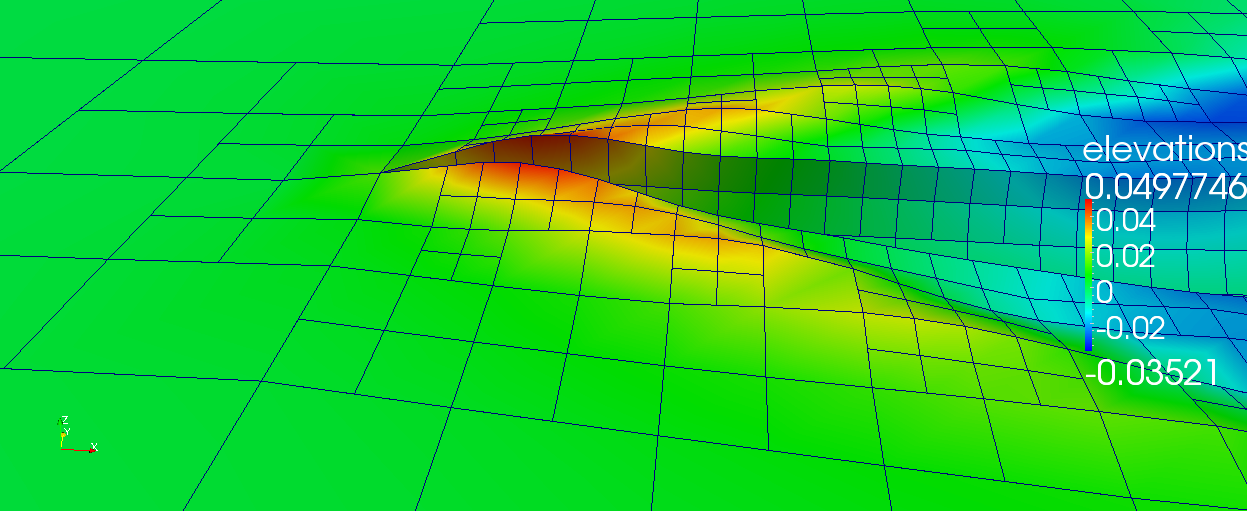}
\end{tabular}
\caption{A mesh refinement step. \label{fig:adaptiveMeshRef}}	
\end{figure}

At fixed intervals in time, the surface gradient of the finite element
approximation of $\phi$ is post-processed. This provides a quantitative
estimate of the cells in which the approximation error may be higher. In
particular, for each cell $K$ of our triangulation we compute the
quantity
\begin{equation}
  \tau_K^2 := \frac{h}{24}\int_{\partial K} \left[ \nablab_s \phi \cdot \nb_{\partial
      K} \right ]^2 \d \gamma,
\end{equation}
where $[ \nablab_s \phi \cdot \nb_{\partial K} ]$ denotes the jump of
the surface gradient of $\phi$ across the edges of the triangulation
element $K$. The vector $\nb_{\partial K}$ is perpendicular to both
the cell normal $\nb$, and to the boundary of the element $K$, and $h$
is the diamater of the cell itself.

Roughly speaking, $\tau_K$ gives an estimate of how well the trial
space is approximating the surface gradient of $\phi$. The higher
these values are, the smaller the cells should become. The estimated error
per cell $\tau_K$ is ordered, and a fixed fraction of the cells with
the highest and lowest error $\tau_K$ are flagged for refinement and
coarsening. The computational grid is then refined, ensuring that any
two neighboring cells differ for at most one refinement level. 

Standard interpolation is used to transfer all finite dimensional
solutions from one grid to another, and a geometrically consistent
mesh modification algorithm is used to collocate the new nodes
coordinates as smoothly as possible
(see~\cite{BonitoNochettoPauletti-2010-a} for a detailed explanation
of this algorithmic strategy). The resulting computational grid is
non conformal. At each hanging node, all the finite dimensional fields
are constrained to be continuous. This results in a set of algebraic
constraints to be applied to all the degrees of freedom associated
with the hanging nodes, which are eliminated from the final system of
equations via a matrix condensation technique.

Most of these algorithmic strategies are based on the ones which were
already implemented in the \texttt{deal.II} finite element
library~\cite{BangerthHartmannKanschat2007} for trial spaces of finite
elements defined in two and three dimensions, and were modified to
allow their use on arbitrary surfaces embedded in three dimensions. An
example of a refinement step is presented in
Figure~\ref{fig:adaptiveMeshRef}.

\begin{figure}[htb!]
  \centering
  \includegraphics[width=.8\textwidth]{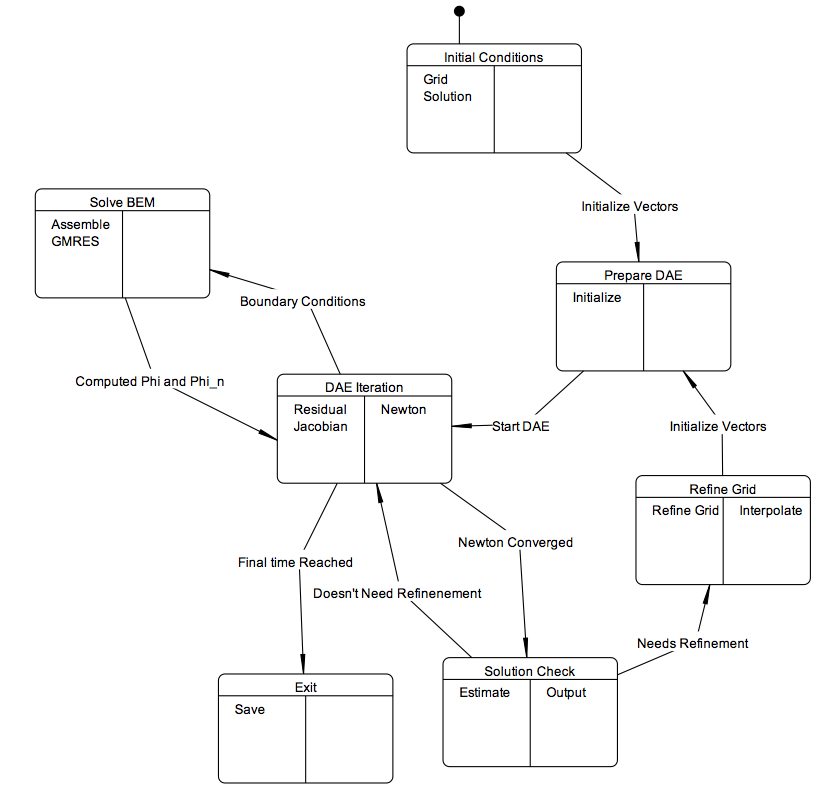}
  \caption{State Diagram}
  \label{fig:state_diagram}
\end{figure}

After each coarsening and refinement step, the system of differential
algebraic equations is restarted with the newly interpolated solution
as initial condition. A state diagram for the entire solution process
is sketched in Figure~\ref{fig:state_diagram}.

\section{Numerical simulations and results}

\begin{figure}[htb]
  \centering
  \includegraphics[width=0.5\textwidth]{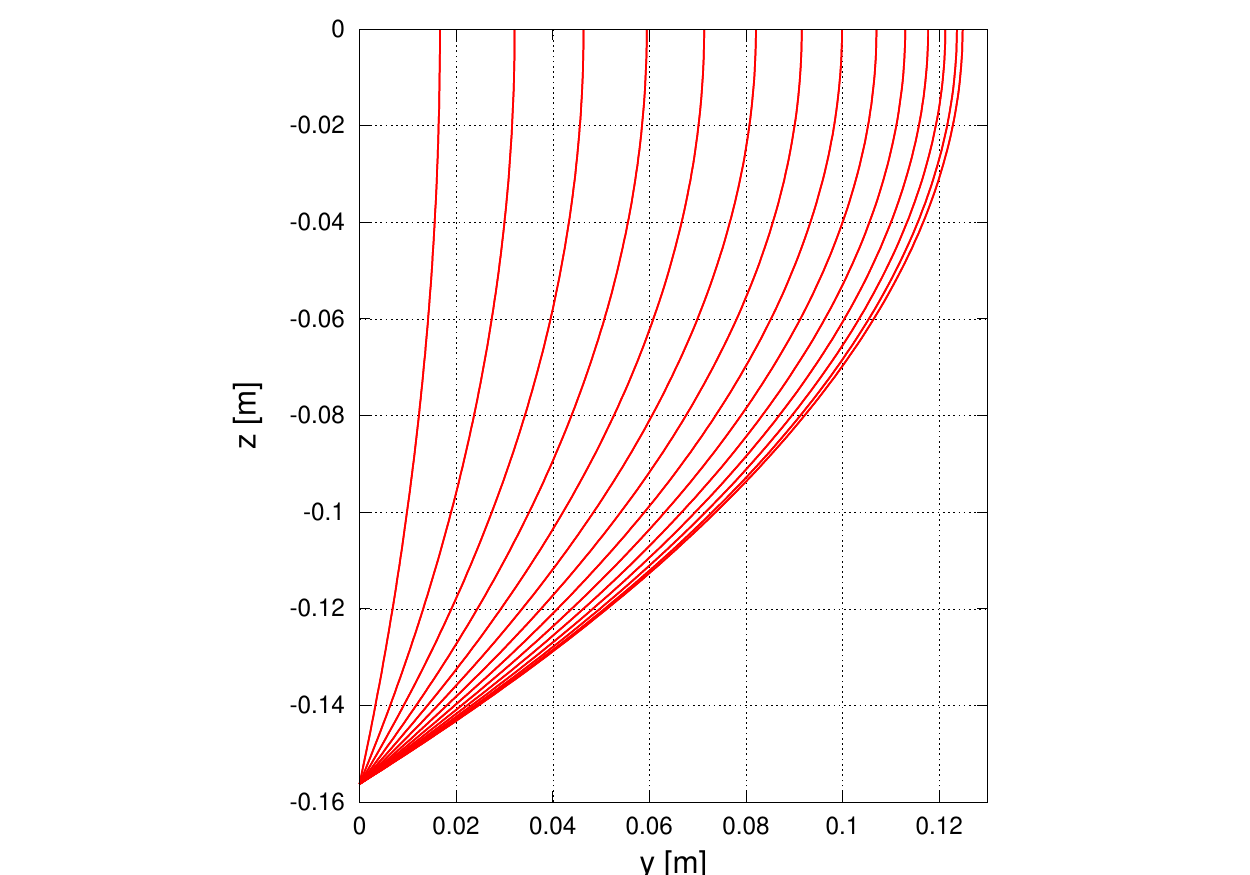}
  \caption{Vertical sections of the Wigley hull used for the simulations, generated by planes
           normal to the longitudinal axis of the hull.}
  \label{fig:wigleyStrips}
\end{figure}

The test case presented in this work is the problem of a Wigley hull
advancing in calm water at speed $\Vb_\infty$ parallel to the
longitudinal axis, with fixed sinkage and trim. In naval engineering,
the Wigley hull is commonly used as a benchmark for the validation of
free-surface flow simulations.  This is mainly due to its simple
shape, defined by the equation

\begin{equation}
y(x,z) = \frac{B}{2}\left[1- \left(\frac{2x}{L}\right)^2\right]\left[1- \left(\frac{z}{T}\right)^2\right].
\end{equation}
In our simulations the boat length, beam and draft values used are respectively  
$L = 2.5\ \text{m}$, $B = 0.25\ \text{m}$, and $T = 0.15625\ \text{m}$. A sketch of the
resulting hull shape is presented in Fig. \ref{fig:wigleyStrips}, which represents a
set of vertical sections of the hull. 

In the numerical simulation setup, the boat is started at rest, and
the its velocity is increased up to the desired speed $\Vb_\infty$
with a linear ramp. The simulation is then carried on until the flow
approaches a steady state solution. The presence of the ramp is not
needed for the stability and convergence of the solution, which are
also obtained imposing an impulsive start of the water past the
hull. Still, inducing slower dynamics in the first seconds of the
simulations, the linear ramp allows for higher time steps and faster
convergence.


To compare the non linear free-surface BEM solutions with the
experimental results presented in~\cite{mccarthy1985}, we considered
six different surge velocities $V_\infty$, corresponding to the Froude
numbers reported in Table \ref{tab:froudeNumbers}. For each of these
Froude numbers, numerical solutions were obtained using a refined and
a coarse mesh, in which the adaptive mesh refinement algorithm was
tuned in order to maintain the cell dimensions over a given minimum
value, and in order to limit the number of degrees of freedom under a
maximum value.

Despite the fact that the BEM solver developed allows for the use of
panels with arbitrary order, the results we present only refer to
iso-parametric bilinear elements. The proposed stabilization mechanism
seems insufficient for higher order quadrilateral elements, especially
when high Froude numbers are considered.  

We are currently investigating alternative stabilization procedures as
well as finer tuning of the SUPG parameters to increase the robustness
of the algorithm when high order panels are used.

\begin{table}[htb]
\centering
\begin{tabular}{||c||c|c|c|c|c|c||}
\hline
\hline
$V_\infty$ &
1.2381$\ \frac{m}{s}$ & 1.3223$\ \frac{m}{s}$ & 1.4312$\ \frac{m}{s}$ &
1.5649$\ \frac{m}{s}$ & 1.7531$\ \frac{m}{s}$ & 2.0205$\ \frac{m}{s}$ \\
\hline
Fr &
0.250 & 0.267 & 0.289 & 0.316 & 0.354 & 0.408 \\
\hline
\hline
\end{tabular}
\caption{Wigley hull surge velocities imposed in each numerical simulation, and the corresponding
         Froude numbers $\text{Fr} = \frac{V_\infty}{\sqrt{gL}}$.
\label{tab:froudeNumbers}}
\end{table}

\begin{figure}[htb!]
  \centering
  
  \subfigure[Fr = 0.250]{\includegraphics[width=\textwidth]{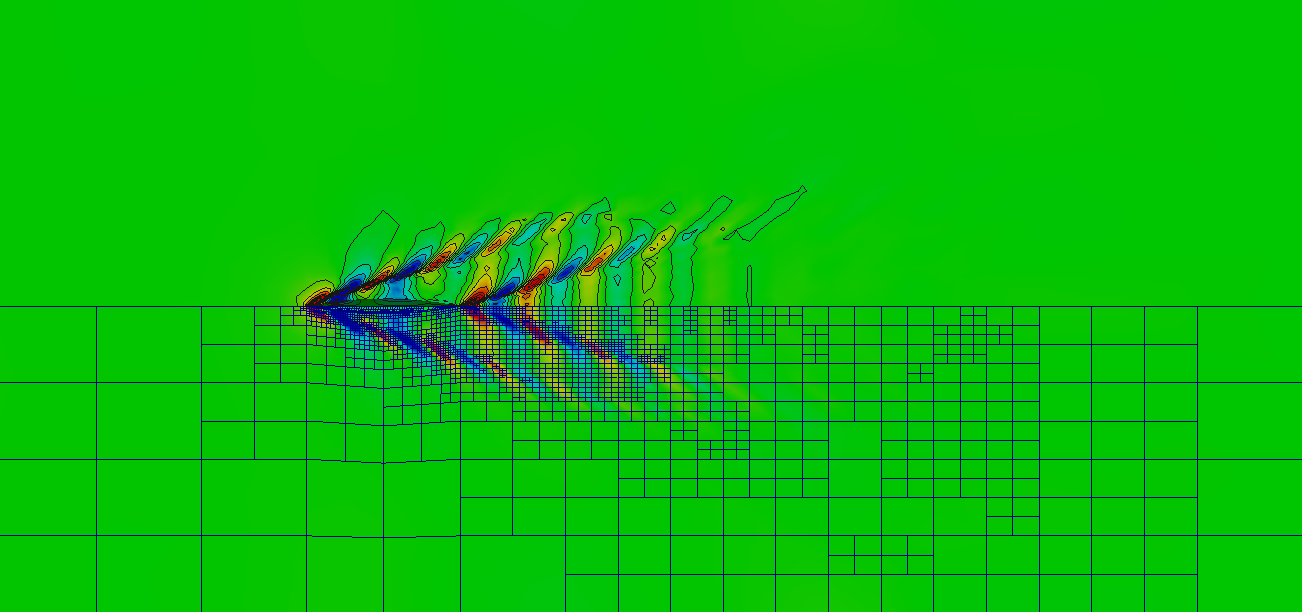}}
  \subfigure[Fr = 0.267]{\includegraphics[width=\textwidth]{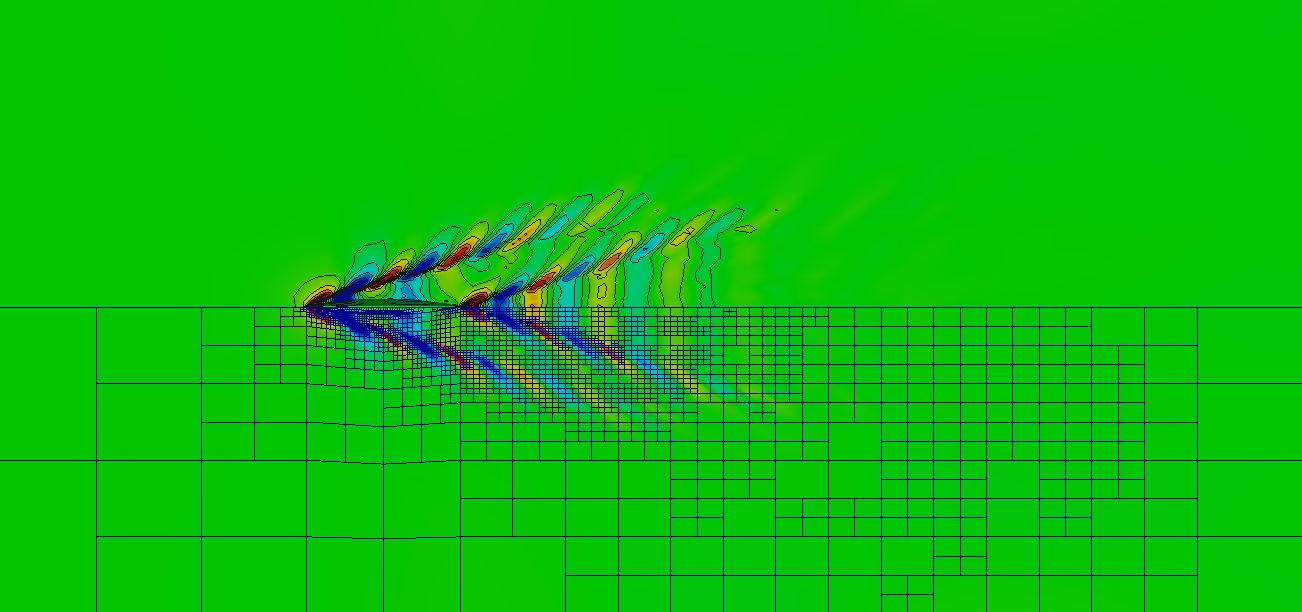}}
  \caption{Mesh refinements and contours (I).}
  \label{fig:mesh-and-contours-1}
\end{figure}

\begin{figure}[htb!]
  \centering
  
  \subfigure[Fr = 0.289]{\includegraphics[width=\textwidth]{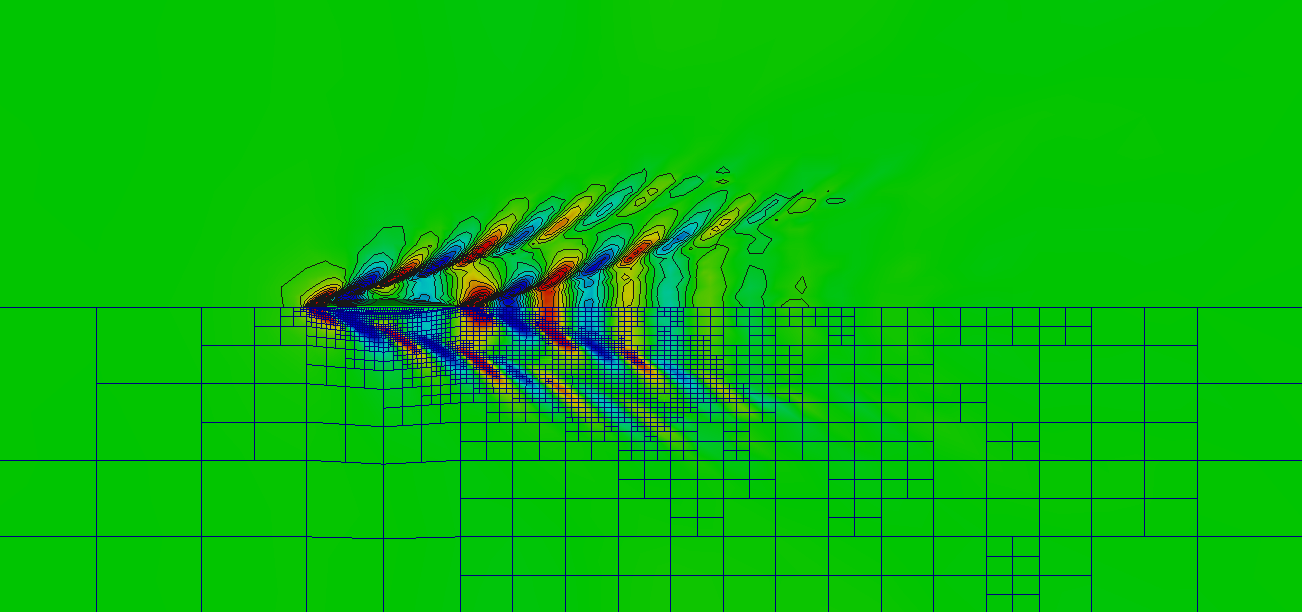}}
  \subfigure[Fr = 0.316]{\includegraphics[width=\textwidth]{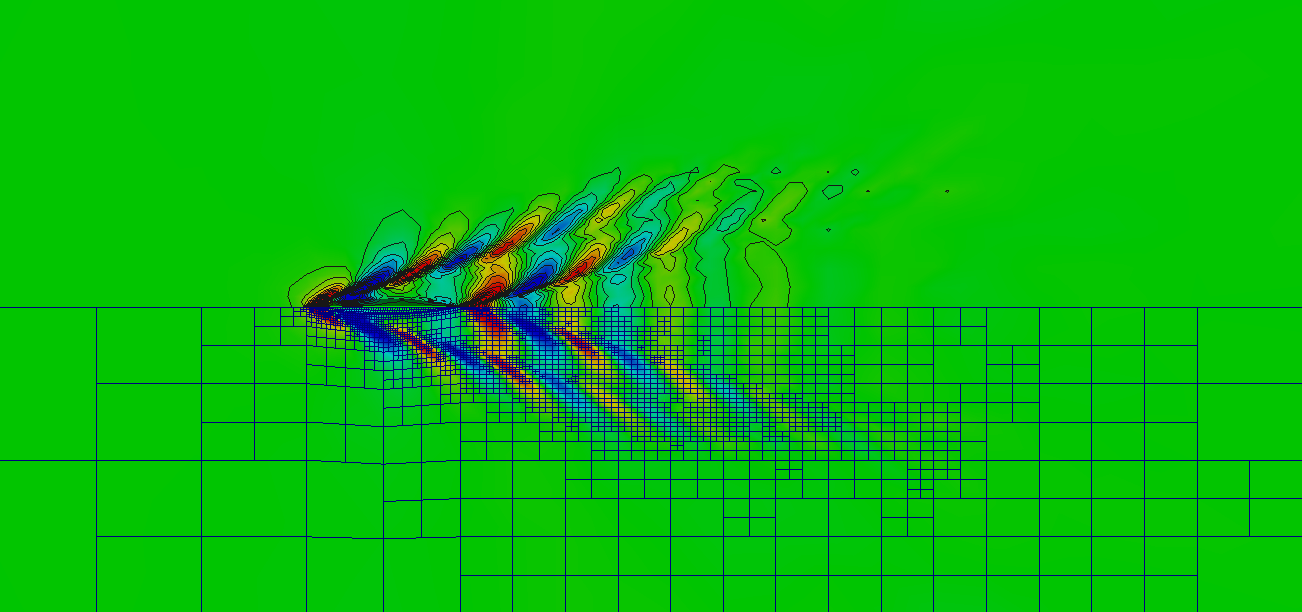}}
 \caption{Mesh refinements and contours (II).}
  \label{fig:mesh-and-contours-2}
\end{figure}

\begin{figure}[htb!]
  \centering

  \subfigure[Fr = 0.354]{\includegraphics[width=\textwidth]{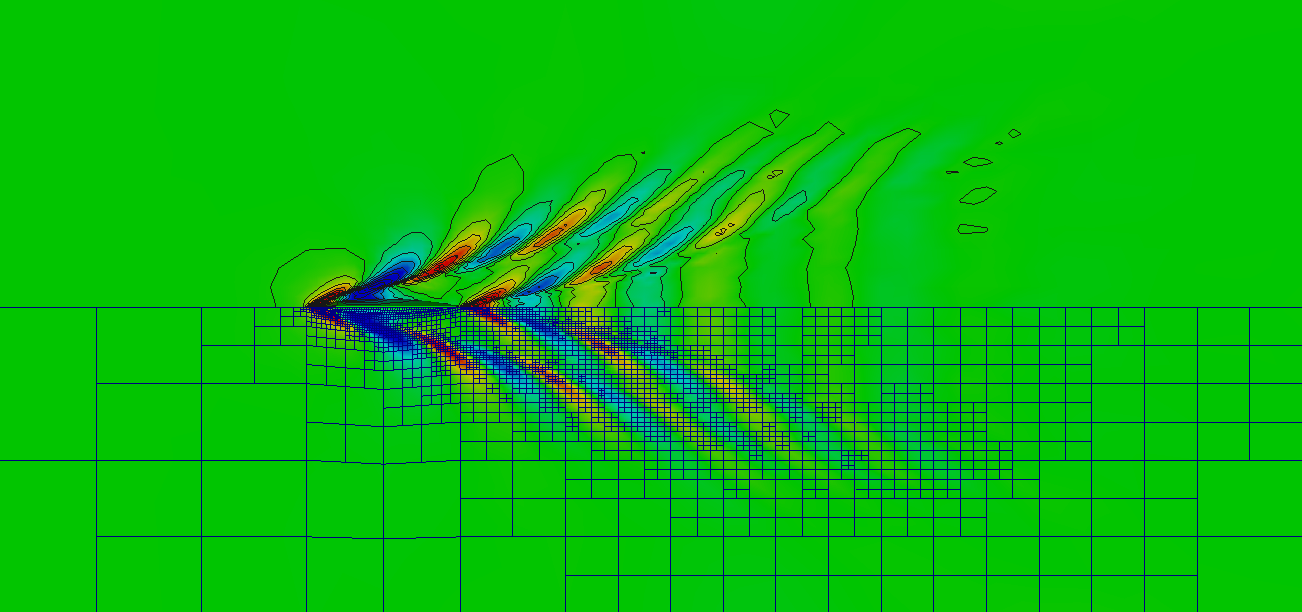}}
  \subfigure[Fr = 0.408]{\includegraphics[width=\textwidth]{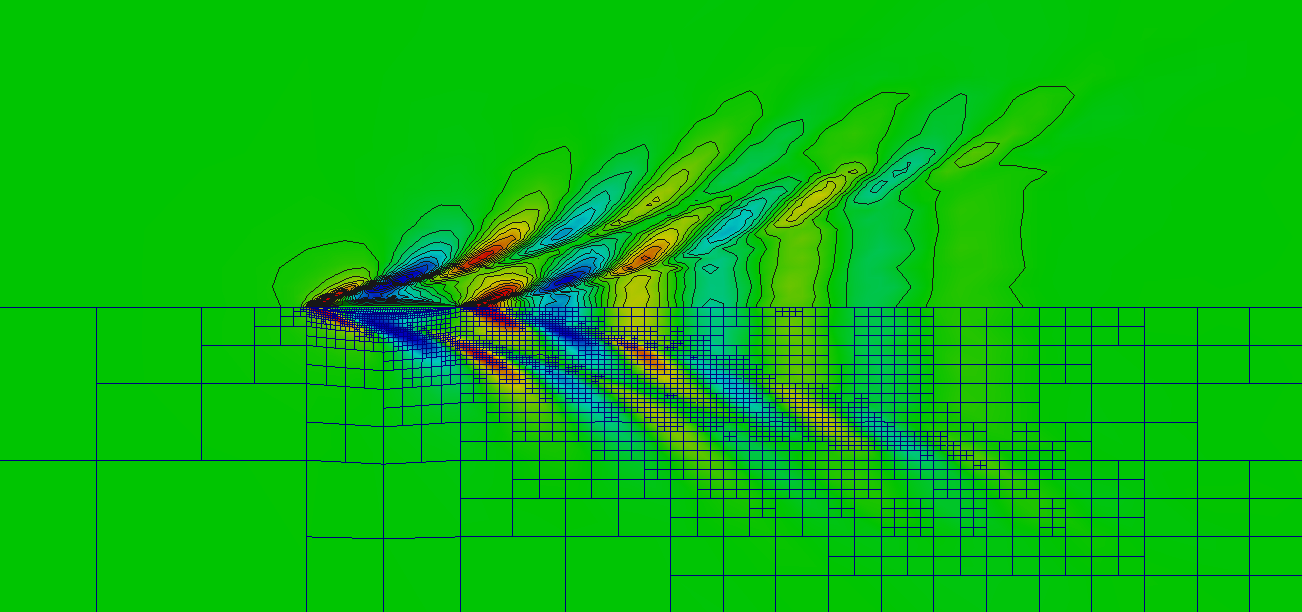}}

  \caption{Mesh refinements and contours (III).}
  \label{fig:mesh-and-contours-3}
\end{figure}

A contour of the wave elevation field for the regime solution obtained
at the various Froude numbers is presented, along with the final mesh, in
Figures~\ref{fig:mesh-and-contours-1}, \ref{fig:mesh-and-contours-2}
and~\ref{fig:mesh-and-contours-3}. The pictures show how the adaptive
mesh refinement leads to an automatic clustering of mesh cells
in regions with highest solution gradients. In this way the numerical
solutions capture in a very accurate manner the
physical characteristics of the wave patterns, requiring a very limited
number of degrees of freedom. The biggest final mesh is in fact only composed
by roughly 6000 nodes, but it allows for a very good reconstruction of
the Kelvin wake, extending for several wavelengths past the
surging hull.

The simulations required 12 hrs (for the coarse meshes) to 48 hrs (for the refined meshes)
to reach the steady state solution, on single SMP nodes of the Arctur cluster of the
Italian/Slovenian interstate cooperation Exact-Lab/Arctur.

\begin{figure}[htb!]
  \centering
\begin{tabular}{c c}
  \includegraphics[width=0.48\textwidth]{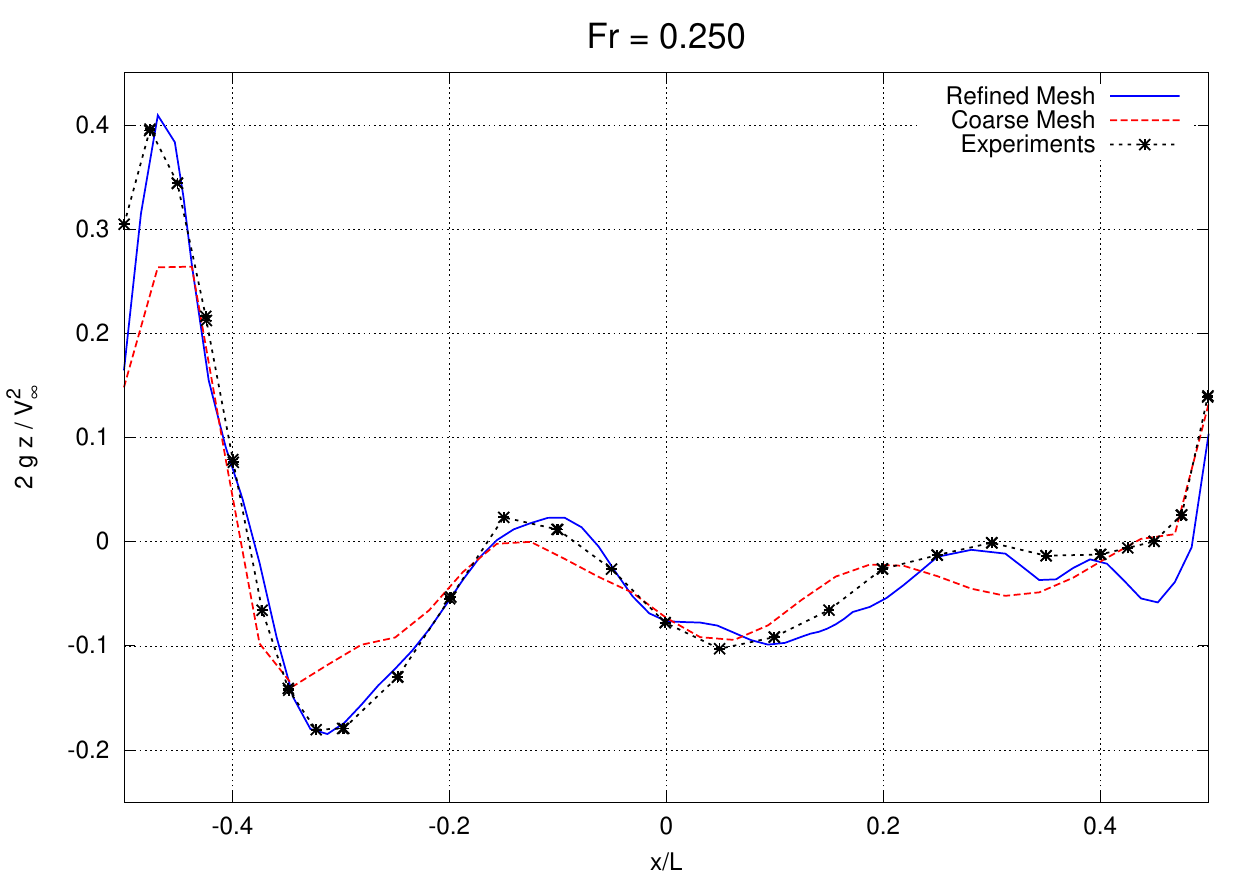}&
  \includegraphics[width=0.48\textwidth]{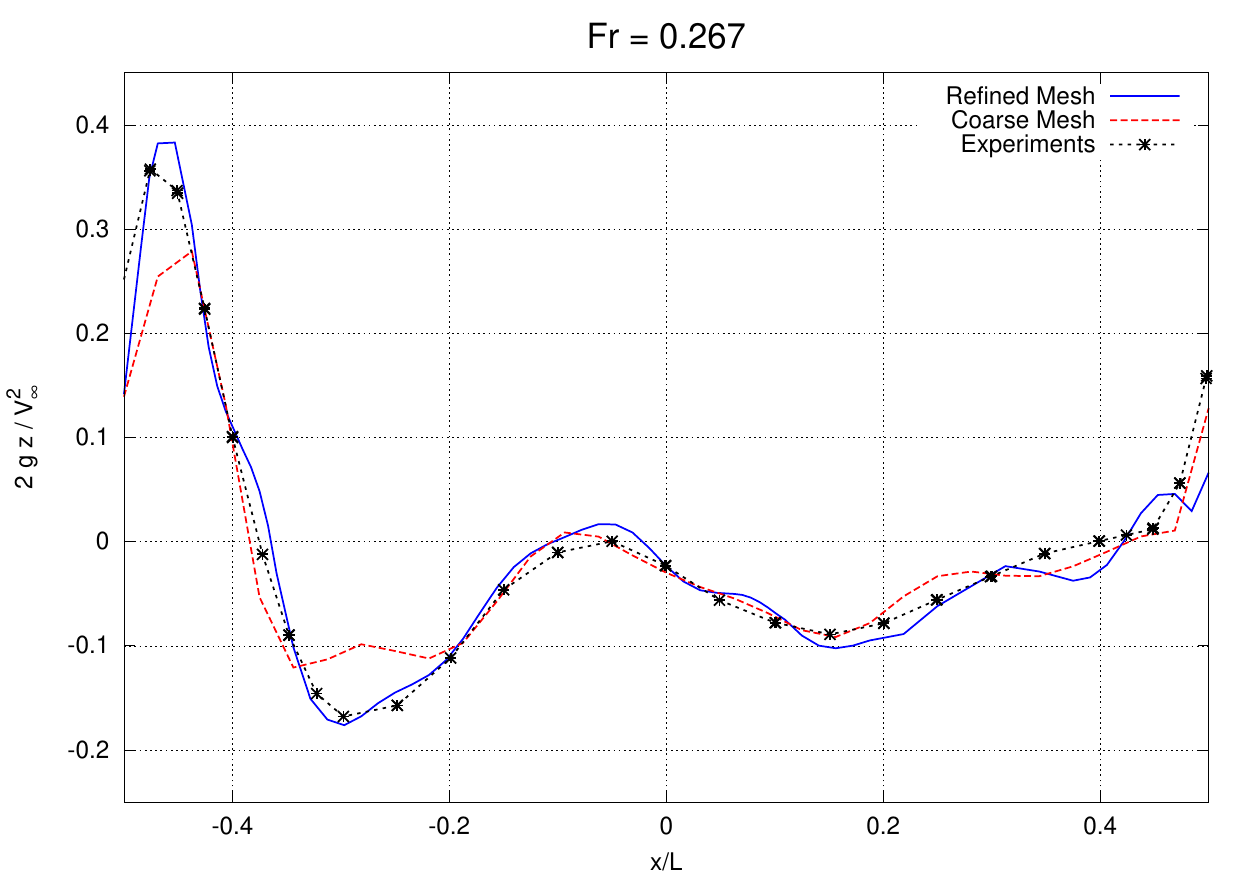}\\
  \includegraphics[width=0.48\textwidth]{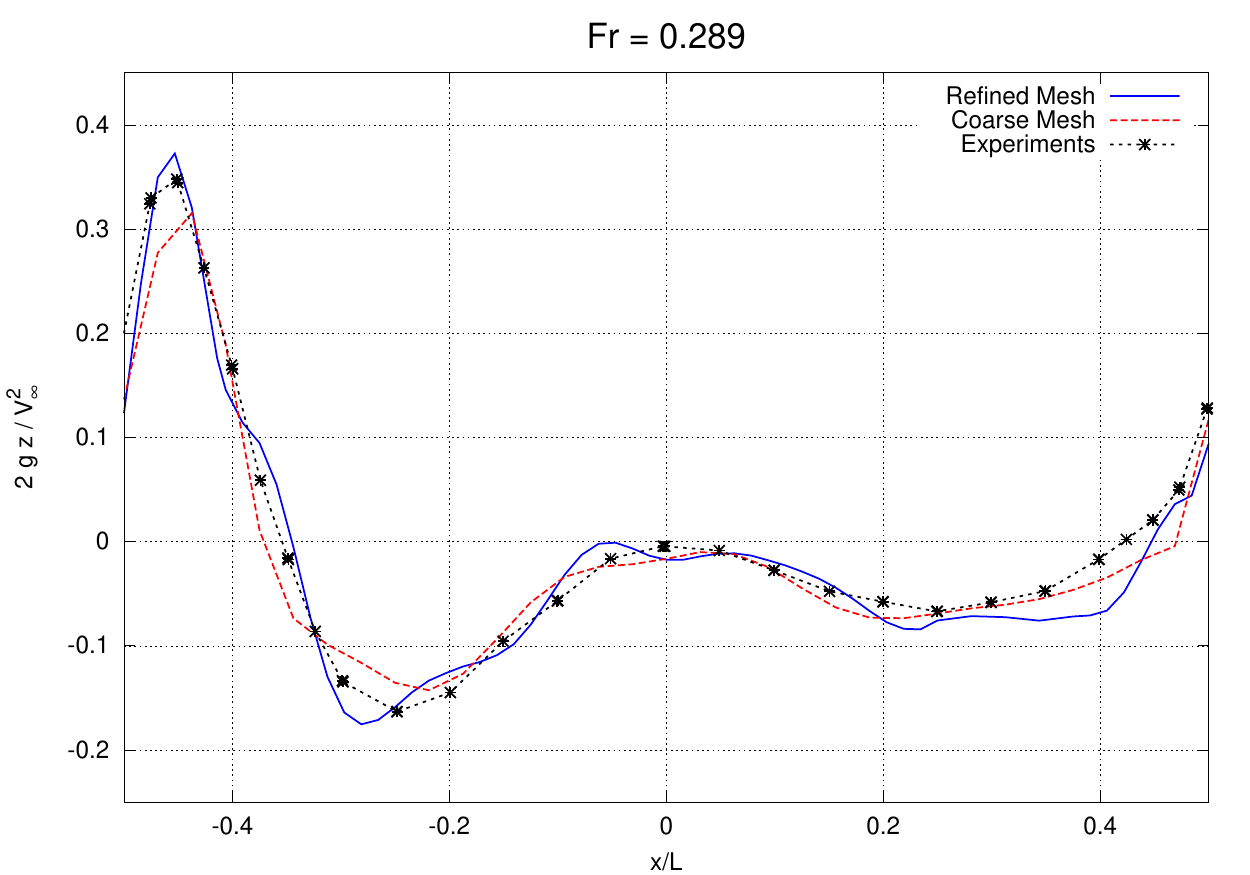}&
  \includegraphics[width=0.48\textwidth]{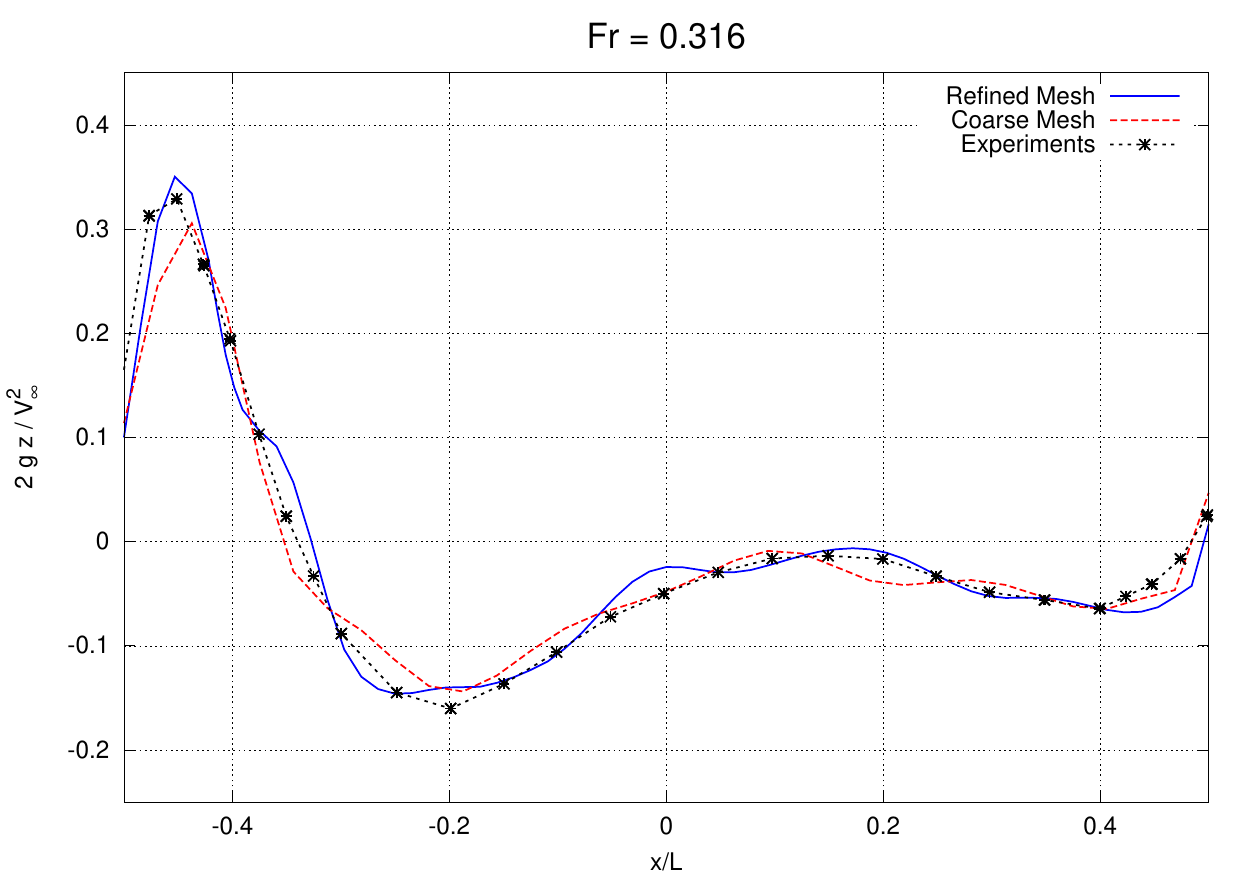}\\
  \includegraphics[width=0.48\textwidth]{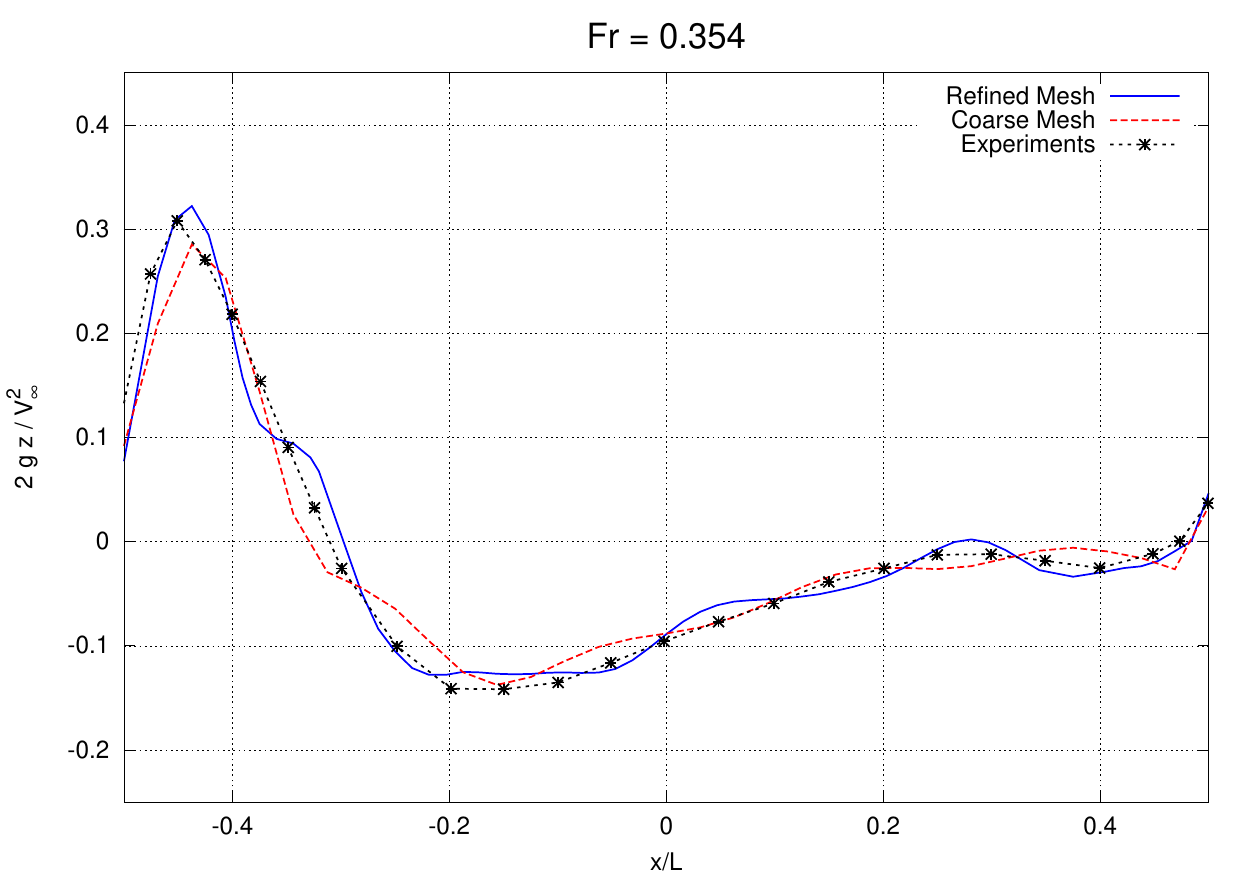}&
  \includegraphics[width=0.48\textwidth]{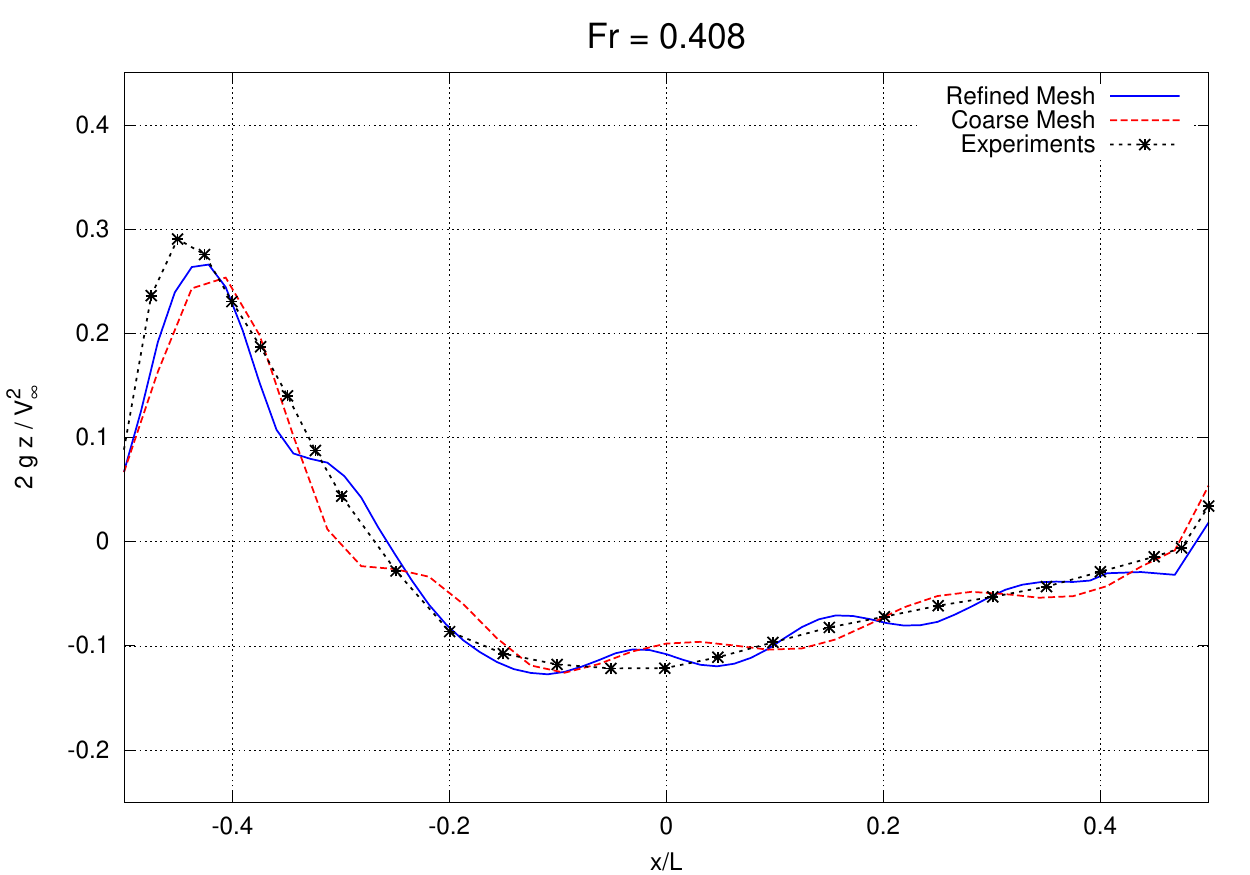}
\end{tabular}

  \caption{Comparison of predicted water profiles with the University of Tokyo
           experimental results ({\color{black}- -$\ast$- -}). Both coarse
           mesh ({\color{red}-- --}) and refined mesh results ({\color{blue}---})
           are shown in the plots.}
  \label{fig:waterLine}
\end{figure}

The wave profiles on the surface of the Wigley hull obtained with the
present method, are compared with the corresponding experimental
results in Fig. \ref{fig:waterLine}. In each plot, the abscissae
represent the dimensionless coordinate $x/L$ along the boat, while the
ordinates are the dimensionless wave elevations $\eta ' =
\frac{2g\eta}{{V_\infty}^2}$. For all the Froude numbers considered,
the method presented seems able to predict qualitatively correct wave
profiles. Moreover, the wave elevation in proximity of the bow of the
boat is reproduced with very good accuracy in all the test cases
considered. On the other hand, in all the numerical curves, we 
observe a small spatial oscillation superimposed to the main wave
profile. The wave length of such oscillation seems proportional to the
local mesh cells size, while the amplitude is slightly higher for
finer meshes. This suggests that better tuning of the SUPG stabilization
parameter may be needed for this kind of boundary value
problems. 

\section{Conclusions}
\label{sec:conclusions}

An accurate and efficient boundary element method for the simulation
of unsteady and fully nonlinear potential waves past surging ships was
developed, implemented and tested. Compared to existing algorithms,
the method presents several innovative features which try to address
some of the most CPU intensive aspects of this kind of computations.

The most innovative idea behind the proposed method is the fact that
the equations are studied on a fixed reference domain, which is
deformed through an arbitrary Lagrangian Eulerian map that keeps track
of the physical shape of the water domain around the ship. Some
aspects of this approach resemble the semi-Lagrangian formulation of
the potential wave equations, but here they are takled using powerful
differential geometry tools, combined with finite element techniques
for arbitrary surfaces.

This reformulation in terms of a fixed reference domain presents
severe stability issues in presence of a forward ship motion, or in
presence of an incident stream velocity. Stabilization is achieved via
a weighted SUPG projection, which allows the use of fully unstructured
meshes, and guarantees an accurate reconstruction of the velocity
fields on the mesh nodes, also when low order finite dimensional
spaces are used for the numerical discretization.

To the best of the authors' knowledge, such formulation has never been
successfully used in ship hydrodynamic problems in presence of a non
zero stream velocity, due to the free-surface instabilities.
 
With respect to existing methods, the combination of the
semi-Lagrangian approach with the SUPG stabilization eliminates the
need for periodic remeshing of the computational domain, and opens up
the possibility to exploit local adaptivity tools, typical of finite
element discretizations.

We exploit these ideas by employing simple a-posteriori error
estimates to refine adaptively the computational mesh. Accurate
results are obtained even when using a very limited number of degrees
of freedom.

Implicit BDF methods with variable order and variable step size are
also employed, which render the final computational tool very
attractive in terms of robustness. 

A direct interface with standard CAD file formats is currently under
development, and our preliminary results indicate that the final tool
could be used to study the unsteady interaction between arbitrary hull
shapes and nonlinear water waves in a robust and automated way.

\section*{Acknowledgments}
\label{sec:aknowledgments}

The research leading to these results has received specific funding
under the ``Young SISSA Scientists' Research Projects'' scheme
2011-2012, promoted by the International School for Advanced Studies
(SISSA), Trieste, Italy.

This work was performed in the context of the project OpenSHIP,
``Simulazioni di fluidodinamica computazionale (CFD) di alta qualita per
le previsioni di prestazioni idrodinamiche del sistema carena-elica in
ambiente OpenSOURCE'', supported by Regione FVG - POR FESR 2007-2013
Obiettivo competitivita regionale e occupazione.

The authors gratefully acknowledge the Italian/Slovenian cross-border
cooperation eXact-Lab/Arctur for the processing power that was made
available on the Arctur cluster.

\appendix

\section{Arbitrary Lagrangian Eulerian Formulation}
\label{sec:ale}

\begin{figure}[htb!]
  \centering
  
  \ifpdf
  \resizebox{1\textwidth}{!}{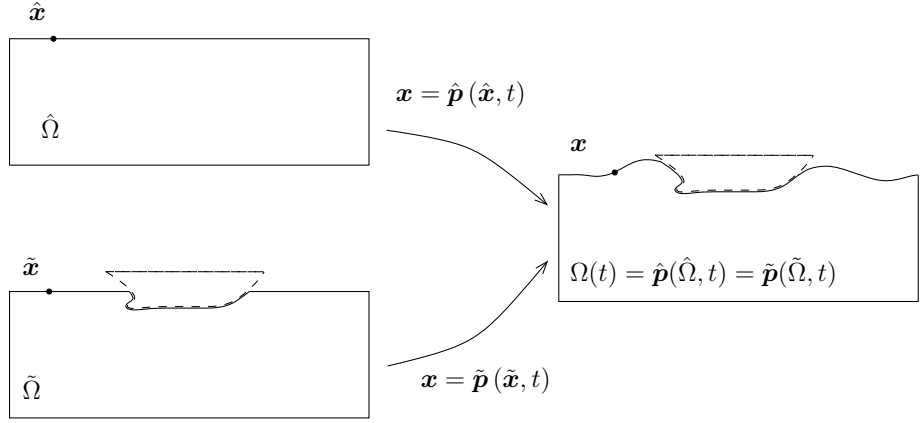}
  \else
  \resizebox{1\textwidth}{!}{\input{figures/ALE.pstex_t}}
  \fi

  \caption{Schematic representation of a Lagrangian motion and of an arbitrary
    Lagrangian Eulerian motion.}
  \label{fig:lag-ale-motions}
\end{figure}



A \emph{motion} is a time parametrized family of invertible maps which
associates to each point $\xlag$ in a reference domain $\Olag$ its
position in space at time $t$:
\begin{equation}
  \label{eq:lagrangian-motion}
  \plag: \Olag\times\Re \mapsto \Re^3, \qquad (\xlag,t) \mapsto
  \xb=\plag(\xlag, t).
\end{equation}
\nomenclature{$\Olag$}{Reference Lagrangian domain}
\nomenclature{$\xlag$}{Coordinates in the reference Lagrangian domain}
\nomenclature{$\plag$}{Deformation of the reference Lagrangian domain}
\nomenclature{$\hat \cdot$}{A material field}
The domain $\Omega(t)$ at the current time can be seen as the image
under the motion $\plag$ of a reference domain $\Olag$, i.e., 
$\plag(\Olag,t)=\Omega(t)$. We will indicate with the $\hat \cdot$
symbol a \emph{material motion}, or a motion in which the points
$\xlag$ label material particles.

If one does not want to follow material particles with the domain
$\Omega(t)$, it is possible to introduce another intermediate motion,
called the \emph{ALE motion}, with which we represent deformations of
the domain $\Omega(t)$:
\begin{equation}
  \label{eq:ale-motion}
  \pale: \Oale\times\Re \mapsto \Re^3, \qquad (\xale,t) \mapsto
  \xb=\pale(\xale, t).
\end{equation}
\nomenclature{$\Oale$}{Arbitrary Lagrangian Eulerian domain}
\nomenclature{$\xale$}{Arbitrary Lagrangian Eulerian coordinates}
\nomenclature{$\pale$}{Arbitrary Lagrangian Eulerian deformation}
\nomenclature{$\tilde \cdot$}{An Arbitrary Lagrangian Eulerian field}
These motions can be rather arbitrary, as long as the \emph{shape} of
the domain $\Omega(t)$ is preserved by the motion itself, i.e.,
$\pale(\Oale,t) = \Omega(t)$. The points labelled with $\xale$ do not,
in general, represent material particles. See
Figure~\ref{fig:lag-ale-motions} for a schematic representation of a
Lagrangian motion and of an ALE motion. 

Given a Lagrangian field $\qlag:\Olag\times\Re\mapsto\Re$, its
Eulerian representation is 
\begin{equation}
  q(\xb,t) := \qlag(\plag^{-1}(\xb,t), t), \qquad \forall \xb \in \Omega(t),
\end{equation}
while, given an Eulerian field $q:\Omega(t)\times \Re\mapsto \Re$, its
Lagrangian representation would be
\begin{equation}
  \qlag(\xlag,t) := q(\plag(\xlag,t), t), \qquad \forall \xlag \in \Olag.
\end{equation}
A similar structure is valid for ALE fields:
\begin{alignat}{2}
  q(\xb,t) &:= \qale(\pale^{-1}(\xb,t), t), \qquad &\forall \xb \in
  \Omega(t),\\
  \qale(\xale,t) &:= q(\pale(\xale,t), t), &\forall \xale \in \Oale.
\end{alignat}

The fluid particle velocity $\vb$ which appears in
Problem~(\ref{eq:incompressible-euler}) is the Eulerian
representation of the particles velocity
\begin{equation}
  \vb(\plag(\xlag,t), t) = \hat\vb(\xlag,t) := \frac{\Pa \plag(\xlag,t)}{\Pa t}.
\end{equation}
In a similar way, we define the Eulerian representation of the
\emph{domain velocity}, or \emph{ALE velocity} the field $\wb$ such
that
\begin{equation}
  \wb(\pale(\xale,t), t) = \hat\wb(\xale,t) := \frac{\Pa \pale(\xale,t)}{\Pa t}.
\end{equation}
\nomenclature{$\wb$}{Eulerian representation of the domain velocity}

Time variations of physical quantities associated with material
particles are computed at fixed $\xlag$, generating the usual material
derivative
\begin{equation}
  \frac{D q(\xb,t)}{Dt} := \left.\frac{\partial q(\plag(\xlag,t),t)}{\partial t}
  \right|_{\xlag=\plag^{-1}(\xb,t)} = \frac{\Pa q(\xb,t)}{\Pa t}+\vb \cdot
  \nablab q(\xb,t).
\end{equation}
\nomenclature{$D/Dt$}{Material derivative}

In a similar fashion, if we compute quantities at fixed ALE point
$\xale$, we obtain the ALE time derivative, wich we will denote as
\begin{equation}
  \frac{\delta q(\xb,t)}{\delta t} := \left.\frac{\partial q(\pale(\xale,t),t)}{\partial t}
  \right|_{\xale=\pale^{-1}(\xb,t)} = \frac{\Pa q(\xb,t)}{\Pa t}+\wb \cdot
  \nablab q(\xb,t).
\end{equation}
\nomenclature{$\delta/\delta t$}{Arbitrary Lagrangian Eulerian derivative}

The ALE deformation allows one to define the equations of motions in
Problem~(\ref{eq:incompressible-euler}) in terms of the fixed ALE
reference domain $\Oale$, while still solving the same physical
problem. On the free-surface part of the boundary, the ALE motion
needs to follow the evolution of the fluid particles in order to
maintain the correct shape of the domain $\Omega(t)$, in particular
the minimum requirement for the ALE motion on the free-surface is
given by the \emph{free-surface kinematic boundary condition}
\begin{equation}
  \label{eq:freesurface-kinematic-condition}
  \wb\cdot\nb = \vb\cdot\nb \qquad \text{ on }\Gamma^w(t),
\end{equation}
\nomenclature{$\nb$}{Unit outer normal vector}
which complements the dynamic boundary
condition~(\ref{eq:freesurface-condition}), and provides an evolution
equation for the normal part of the ALE motion on the boundary
$\Gamma^w(t)$. In terms of the ALE motion, the above condition reads
\begin{equation}
  \label{eq:freesurface-kinematic-condition-ale}
  \frac{\Pa\pale}{\Pa t}(\xale, t) \cdot \nb = \vb(\pale(\xale, t),
  t)\cdot\nb \qquad \text{ on }\Gamma^w(t).
\end{equation}

Equations~(\ref{eq:freesurface-condition})
and~(\ref{eq:freesurface-kinematic-condition}) are usually referred to
as the free-surface dynamic and kinematic boundary conditions, since
they represent the physical condition applied to the free-surface
(equilibrium of the pressure on the water surface) and its evolution
equation (the shape of the free-surface follows the velocity field of
the flow).

\section{Numerical beach}
\label{sec:numerical-beach}

A draw back of using an homogeneous Neumann boundary conditions for
the vertical far field boundary condition is that it reflects energy
back in the computational domain.

We use an \emph{absorbing beach} technique (see, for example,
\cite{CaoBeckSchultz-1993}), in which we add an artificial damping
region away from the boat, used to absorb the wave energy. A damping
term can be seen as an additional pressure $P$ acting on the free
surface. In such case, Bernoulli equation on the free-surface becomes

\begin{equation}
\frac{\Pa\Phi}{\Pa t}+gz- \ab_{\stream}\cdot\xb +\frac{1}{2}|\nablab\Phi|^2 + \frac{P}{\rho} = 0
\qquad\qquad \text{ on } \ \ \Gamma^w(t),
\label{eq:bernoulliDamp}
\end{equation}
and one can show that the resulting rate of energy absorption is
\nomenclature{$P$}{Damping pressure for numerical beach}
\begin{equation}
\frac{dE_f}{dt} = \int_{\Gamma^w}P\phin\d\Gamma.
\label{eq:dampRate}
\end{equation}
A natural way to construct the damping pressure $P$ is then to use a
term which is proportional to the potential normal derivative
$\phin$, which grants a positive energy absorption at
all times.

The dynamical free-surface boundary condition modified to account for the damping
term reads 
\begin{equation}
\frac{\Pa\Phi}{\Pa t}+gz- \ab_{\stream}\cdot\xb +\frac{1}{2}|\nablab\Phi|^2 -
\mu\phin = 0
\qquad\qquad \text{ on } \ \ \Gamma^w(t),
\label{eq:fsDynamicBeckDamped}
\end{equation}
where
\begin{equation}
\mu = \left(\frac{\max\left(0,x-x_d\right)}{L_d}\right)^2,
\end{equation}
and $x_d$ is the $x$ coordinate value in which the artificial
damping starts to act, while $L_d$ is the length of the absorbing
beach, as in Figure~\ref{fig:wholeDomain}.

\nomenclature{$x_d$}{$x$-coordinate of the starting point of the
  numerical beach}
\nomenclature{$\nu$}{Damping coefficient of the numerical beach}

\section{Iso-parametric discretization}
\label{sec:iso-param}

We approximate the geometry of the domain boundary by means of
arbitrary order quadrilateral panels. The edges of each panel are
defined by polynomial functions, while their internal surface
is described by polynomial tensor products.

\begin{figure}[htb!]        
\centerline{
  \ifpdf
  \resizebox{1\textwidth}{!}{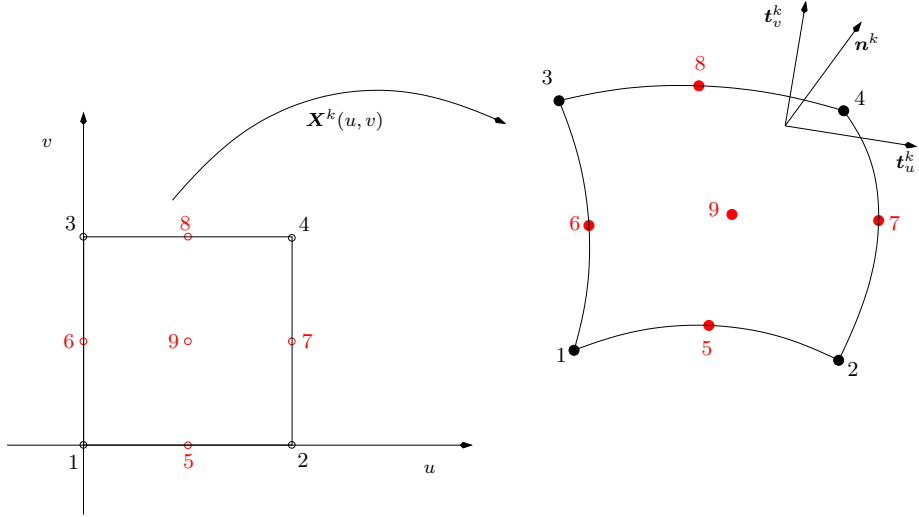}
  \else
  \resizebox{1\textwidth}{!}{\input{figures/pannello.pstex_t}}
  \fi
}
\caption{A linear panel and the reference element. The red dots represent
         the additional degrees of freedom of a quadratic panel. \label{fig:panel}}	
\end{figure}  

In particular we employ the Lagrangian shape functions $N_l(u,v) \quad
l = 1,\dots,N_L$ on the reference panel (Fig. \ref{fig:panel}, on
the left). Thus, the local parametrization of the $k$-th
panel reads
\begin{equation}
  \label{eq:local-gamma-h}
  \begin{split}
    \xale_k(u,v) &:= \sum_{l=1}^{N_L} \xale^{k_l} N_l(u,v) \qquad \qquad
    u,v\in[0,1]^2\\
    \xb_k(u,v,t) &:= \sum_{l=1}^{N_L} \xb^{k_l}(t) N_l(u,v) \qquad \qquad
    u,v\in[0,1]^2.
  \end{split}
\end{equation}
The parametrization weights in Eq. (\ref{eq:local-gamma-h}) are the positions
of the nodes in the reference domain $\tilde\Gamma_h$, or in the current
domain $\Gamma_h(t)$. $k_l$ denotes the \emph{local to global} numbering index
which identifies the $N_L$ basis functions $\varphi^{k_l}$ which are
different from zero on the $k$-th panel.

The global basis functions $\varphi^i(\xale)$ can be identified
and evaluated on each panel $K$ via their local parametrization as
\begin{equation}
  \label{eq:global-basis}
  \varphi^i_k(u,v) := \varphi^i(\xale_k(u,v)) =  \sum_{l=1}^{N_L}\delta_{i\,k_l}N_l(u,v),\qquad 
  \delta_{ij} =
  \begin{cases}
    1 & \text{ if } i = j\\
    0 & \text{ otherwise.}
  \end{cases}
\end{equation}

In this framework, the local representation of $\phi(\xb^k(u,v,t),t)$ and of
its normal derivative on the $k$-th panel are
\begin{equation}
\phi_k(u,v,t) = \sum_{l=1}^{N_L} \phi^{k_l}(t) N_l(u,v) \qquad
\phin{}_k(u,v,t) =
\sum_{l=1}^{N} \phin ^{k_l}(t)
N_l(u,v),
\end{equation}
where $\phi^{k_l}, \phin^{k_l}
\,\,l=1,\dots,N_L$ are the nodal values of the potential and of its
normal derivative in panel $k$.  
On each point of the panel, it is possible to compute two vectors
tangential to $\Gamma_h(t)$ as
\begin{equation}
\tb_u^k(u,v,t) = \sum_{l=1}^{N_L} \xb^{k_l}(t) \frac{\Pa N_l}{\Pa u}(u,v)
\qquad
\tb_v^k(u,v,t) = \sum_{l=1}^{N_L} \xb^{k_l}(t) \frac{\Pa N_l}{\Pa v}(u,v),
\end{equation}
%
from which the external normal unit vector $\nb$ is obtained as
\begin{equation} \label{eq:normalComp}
\nb_{vect}^k(u,v,t) = \tb_u^k(u,v,t) \times \tb_v^k(u,v,t),
\qquad
\nb^k(u,v,t) = \frac{\nb_{vect}^k(u,v,t)}
                    {\left|\nb_{vect}^k(u,v,t)\right|}.
\end{equation}

\nomenclature{$\tb_{u}^k$}{Tangential vector to the $k$-th panel in
  the $u$ direction} 
\nomenclature{$\tb_{v}^k$}{Tangential vector to the $k$-th panel in
  the $v$ direction} 
\nomenclature{$\nb_{vect}^k$}{Normal vector to the $k$-th panel}
\nomenclature{$\nb^k$}{Unit normal vector to the $k$-th panel}

The same can be done for vectors tangential and normal to
$\tilde\Gamma$, by simply replacing $\xb^{k_l}(t)$ with $\xale^{k_l}$
in the definitions above. We will denote those vectors with
$\tilde\tb_u^k(u,v)$, $\tilde\tb_v^k(u,v)$, $\tilde\nb_{vect}^k(u,v)$
and $\tilde\nb^k(u,v)$, respectively.

Integrals on a panel $K$ (or $\tilde K$), can be computed in the
reference domain $[0,1]^2$, by observing that $\d \Gamma = J^k(u,v,t)
\d u\d v$, where $J^k(u,v,t):=|\nb_{vect}^k(u,v,t)|$ (or $\d
\tilde\Gamma = \tilde J^k(u,v) \d u\d v$, where $\tilde
J^k(u,v):=|\tilde\nb_{vect}^k(u,v)|$).

We finally introduce the following differential operators

\begin{equation}
  \begin{split}
    D_k(u,v,t) &:= \nablab_{uv} \xb_k(u,v,t)
    \qquad\in\Re^{3\times2}\\
    G_k(u,v,t) &:= D_k(u,v,t)^T
    D_k(u,v,t)\qquad \in\Re^{2\times2}\\
    \tilde D_k(u,v) &:= \nablab_{uv} \xale_k(u,v)
    \qquad\in\Re^{3\times2}\\
    \tilde G_k(u,v) &:= \tilde D_k(u,v)^T
    \tilde D_k(u,v)\qquad \in\Re^{2\times2},
  \end{split}
\end{equation}
where $G_k$ and $\tilde G_k$ are the first fundamental forms in the
element $K$ and $\tilde K$. Making use of such operators, we can express
the local surface gradients of the global basis functions as 

\begin{equation}
  \label{eq:local-surface-gradient-basis-functions}
  \begin{split}
    \left.\nablab_s\varphi^i(\pale^{-1}(\xb,t))\right.|_{\xb=\xb_k(u,v,t)}
    & = D_k(G_k)^{-1} \nabla_{uv} \varphi^i_k(u,v,t) \\
    & =: \nablab_s  \varphi^i_k(u,v,t)\\
    \left.\tilde\nablab_s\varphi^i(\xale)\right.|_{\xale=\xale_k(u,v)}
    & = \tilde D_k(G_k)^{-1} \nabla_{uv} \varphi^i_k(u,v) \\
    & =: \tilde \nablab_s  \varphi^i_k(u,v),
  \end{split}
\end{equation}
which we will indicate with the same symbol as the spatial surface
gradients. The surface gradient of a finite dimensional vector can then
be computed by Equation~(\ref{eq:surface-discrete-gradient}).

\nomenclature{$D_k$}{First fundamental form in the $k$-th panel}
\nomenclature{$\tilde D_k$}{First fundamental form in the $k$-th
  reference panel}

\section{Nomenclature}
\label{sec:nomenclature}

\printnomenclature

\section*{Bibliography}

\bibliographystyle{plain}
\bibliography{References}

\end{document}

%% file: figures/dominio_barca.pdftex_t
\begin{picture}(0,0)%
\includegraphics{figures/dominio_barca.pdf}%
\end{picture}%
\setlength{\unitlength}{2693sp}%
\begingroup\makeatletter\ifx\SetFigFontNFSS\undefined%
\gdef\SetFigFontNFSS#1#2#3#4#5{%
  \reset@font\fontsize{#1}{#2pt}%
  \fontfamily{#3}\fontseries{#4}\fontshape{#5}%
  \selectfont}%
\fi\endgroup%
\begin{picture}(8858,7071)(2303,-7890)
\put(2632,-7382){\makebox(0,0)[lb]{\smash{{\SetFigFontNFSS{8}{9.6}{\familydefault}{\mddefault}{\updefault}{\color[rgb]{0,0,0}$\Gamma^{b}$}%
}}}}
\put(2946,-5754){\makebox(0,0)[lb]{\smash{{\SetFigFontNFSS{8}{9.6}{\familydefault}{\mddefault}{\updefault}{\color[rgb]{0,0,0}$\Vb_{\infty}$}%
}}}}
\put(10356,-3714){\makebox(0,0)[lb]{\smash{{\SetFigFontNFSS{8}{9.6}{\familydefault}{\mddefault}{\updefault}{\color[rgb]{0,0,0}$L_{d}$}%
}}}}
\put(5639,-3134){\makebox(0,0)[lb]{\smash{{\SetFigFontNFSS{8}{9.6}{\familydefault}{\mddefault}{\updefault}{\color[rgb]{0,0,0}$y$}%
}}}}
\put(7469,-3314){\makebox(0,0)[lb]{\smash{{\SetFigFontNFSS{8}{9.6}{\familydefault}{\mddefault}{\updefault}{\color[rgb]{0,0,0}$x$}%
}}}}
\put(6569,-2789){\makebox(0,0)[lb]{\smash{{\SetFigFontNFSS{8}{9.6}{\familydefault}{\mddefault}{\updefault}{\color[rgb]{0,0,0}$z$}%
}}}}
\put(2543,-1871){\makebox(0,0)[lb]{\smash{{\SetFigFontNFSS{8}{9.6}{\familydefault}{\mddefault}{\updefault}{\color[rgb]{0,0,0}$\Gamma^{w}$}%
}}}}
\put(9247,-7037){\makebox(0,0)[lb]{\smash{{\SetFigFontNFSS{8}{9.6}{\familydefault}{\mddefault}{\updefault}{\color[rgb]{0,0,0}$\Gamma^{ff}$}%
}}}}
\put(4561,-1501){\makebox(0,0)[lb]{\smash{{\SetFigFontNFSS{8}{9.6}{\familydefault}{\mddefault}{\updefault}{\color[rgb]{0,0,0}$\Gamma^h$}%
}}}}
\put(8796,-4434){\makebox(0,0)[lb]{\smash{{\SetFigFontNFSS{8}{9.6}{\familydefault}{\mddefault}{\updefault}{\color[rgb]{0,0,0}$x_{d}$}%
}}}}
\end{picture}%

%% file: figures/dominio_barca.pstex_t
\begin{picture}(0,0)%
\includegraphics{figures/dominio_barca.eps}%
\end{picture}%
\setlength{\unitlength}{2693sp}%
\begingroup\makeatletter\ifx\SetFigFontNFSS\undefined%
\gdef\SetFigFontNFSS#1#2#3#4#5{%
  \reset@font\fontsize{#1}{#2pt}%
  \fontfamily{#3}\fontseries{#4}\fontshape{#5}%
  \selectfont}%
\fi\endgroup%
\begin{picture}(8858,7071)(2303,-7890)
\put(2632,-7382){\makebox(0,0)[lb]{\smash{{\SetFigFontNFSS{8}{9.6}{\familydefault}{\mddefault}{\updefault}{\color[rgb]{0,0,0}$\Gamma^{b}$}%
}}}}
\put(2946,-5754){\makebox(0,0)[lb]{\smash{{\SetFigFontNFSS{8}{9.6}{\familydefault}{\mddefault}{\updefault}{\color[rgb]{0,0,0}$\Vb_{\infty}$}%
}}}}
\put(10356,-3714){\makebox(0,0)[lb]{\smash{{\SetFigFontNFSS{8}{9.6}{\familydefault}{\mddefault}{\updefault}{\color[rgb]{0,0,0}$L_{d}$}%
}}}}
\put(5639,-3134){\makebox(0,0)[lb]{\smash{{\SetFigFontNFSS{8}{9.6}{\familydefault}{\mddefault}{\updefault}{\color[rgb]{0,0,0}$y$}%
}}}}
\put(7469,-3314){\makebox(0,0)[lb]{\smash{{\SetFigFontNFSS{8}{9.6}{\familydefault}{\mddefault}{\updefault}{\color[rgb]{0,0,0}$x$}%
}}}}
\put(6569,-2789){\makebox(0,0)[lb]{\smash{{\SetFigFontNFSS{8}{9.6}{\familydefault}{\mddefault}{\updefault}{\color[rgb]{0,0,0}$z$}%
}}}}
\put(2543,-1871){\makebox(0,0)[lb]{\smash{{\SetFigFontNFSS{8}{9.6}{\familydefault}{\mddefault}{\updefault}{\color[rgb]{0,0,0}$\Gamma^{w}$}%
}}}}
\put(9247,-7037){\makebox(0,0)[lb]{\smash{{\SetFigFontNFSS{8}{9.6}{\familydefault}{\mddefault}{\updefault}{\color[rgb]{0,0,0}$\Gamma^{ff}$}%
}}}}
\put(4561,-1501){\makebox(0,0)[lb]{\smash{{\SetFigFontNFSS{8}{9.6}{\familydefault}{\mddefault}{\updefault}{\color[rgb]{0,0,0}$\Gamma^h$}%
}}}}
\put(8796,-4434){\makebox(0,0)[lb]{\smash{{\SetFigFontNFSS{8}{9.6}{\familydefault}{\mddefault}{\updefault}{\color[rgb]{0,0,0}$x_{d}$}%
}}}}
\end{picture}%

%% file: figures/edgesTreat.pdftex_t
\begin{picture}(0,0)%
\includegraphics{figures/edgesTreat.pdf}%
\end{picture}%
\setlength{\unitlength}{2565sp}%
\begingroup\makeatletter\ifx\SetFigFont\undefined%
\gdef\SetFigFont#1#2#3#4#5{%
  \reset@font\fontsize{#1}{#2pt}%
  \fontfamily{#3}\fontseries{#4}\fontshape{#5}%
  \selectfont}%
\fi\endgroup%
\begin{picture}(4360,3852)(2611,-5680)
\put(4666,-2011){\makebox(0,0)[lb]{\smash{{\SetFigFont{8}{9.6}{\familydefault}{\mddefault}{\updefault}{\color[rgb]{0,0,0}$\nb_1$}%
}}}}
\put(2626,-3481){\makebox(0,0)[lb]{\smash{{\SetFigFont{8}{9.6}{\familydefault}{\mddefault}{\updefault}{\color[rgb]{0,0,0}$\nb_2$}%
}}}}
\end{picture}%

%% file: figures/edgesTreat.pstex_t
\begin{picture}(0,0)%
\includegraphics{figures/edgesTreat.eps}%
\end{picture}%
\setlength{\unitlength}{2565sp}%
\begingroup\makeatletter\ifx\SetFigFont\undefined%
\gdef\SetFigFont#1#2#3#4#5{%
  \reset@font\fontsize{#1}{#2pt}%
  \fontfamily{#3}\fontseries{#4}\fontshape{#5}%
  \selectfont}%
\fi\endgroup%
\begin{picture}(4360,3840)(2611,-5680)
\put(2626,-3481){\makebox(0,0)[lb]{\smash{{\SetFigFont{8}{9.6}{\familydefault}{\mddefault}{\updefault}{\color[rgb]{0,0,0}$\nb_2$}%
}}}}
\put(4666,-2011){\makebox(0,0)[lb]{\smash{{\SetFigFont{8}{9.6}{\familydefault}{\mddefault}{\updefault}{\color[rgb]{0,0,0}$\nb_1$}%
}}}}
\end{picture}%

%% file: figures/ALE.pdftex_t
\begin{picture}(0,0)%
\includegraphics{figures/ALE.pdf}%
\end{picture}%
\setlength{\unitlength}{4144sp}%
\begingroup\makeatletter\ifx\SetFigFont\undefined%
\gdef\SetFigFont#1#2#3#4#5{%
  \reset@font\fontsize{#1}{#2pt}%
  \fontfamily{#3}\fontseries{#4}\fontshape{#5}%
  \selectfont}%
\fi\endgroup%
\begin{picture}(6597,3075)(2596,-5023)
\put(6661,-4021){\makebox(0,0)[lb]{\smash{{\SetFigFont{12}{14.4}{\familydefault}{\mddefault}{\updefault}{\color[rgb]{0,0,0}$\Omega(t)=\plag(\Olag,t)=\pale(\Oale,t)$}%
}}}}
\put(2836,-2986){\makebox(0,0)[lb]{\smash{{\SetFigFont{12}{14.4}{\familydefault}{\mddefault}{\updefault}{\color[rgb]{0,0,0}$\Olag$}%
}}}}
\put(2701,-4876){\makebox(0,0)[lb]{\smash{{\SetFigFont{12}{14.4}{\familydefault}{\mddefault}{\updefault}{\color[rgb]{0,0,0}$\Oale$}%
}}}}
\put(2701,-3976){\makebox(0,0)[lb]{\smash{{\SetFigFont{12}{14.4}{\familydefault}{\mddefault}{\updefault}{\color[rgb]{0,0,0}$\xale$}%
}}}}
\put(2746,-2131){\makebox(0,0)[lb]{\smash{{\SetFigFont{12}{14.4}{\familydefault}{\mddefault}{\updefault}{\color[rgb]{0,0,0}$\xlag$}%
}}}}
\put(5581,-4786){\makebox(0,0)[lb]{\smash{{\SetFigFont{12}{14.4}{\familydefault}{\mddefault}{\updefault}{\color[rgb]{0,0,0}$\xb=\pale \left(\xale,t \right)$}%
}}}}
\put(5401,-2716){\makebox(0,0)[lb]{\smash{{\SetFigFont{12}{14.4}{\familydefault}{\mddefault}{\updefault}{\color[rgb]{0,0,0}$\xb=\plag\left(\xlag,t\right)$}%
}}}}
\put(6661,-3076){\makebox(0,0)[lb]{\smash{{\SetFigFont{12}{14.4}{\familydefault}{\mddefault}{\updefault}{\color[rgb]{0,0,0}$\xb$}%
}}}}
\end{picture}%

%% file: figures/ALE.pstex_t
\begin{picture}(0,0)%
\includegraphics{figures/ALE.eps}%
\end{picture}%
\setlength{\unitlength}{4144sp}%
\begingroup\makeatletter\ifx\SetFigFont\undefined%
\gdef\SetFigFont#1#2#3#4#5{%
  \reset@font\fontsize{#1}{#2pt}%
  \fontfamily{#3}\fontseries{#4}\fontshape{#5}%
  \selectfont}%
\fi\endgroup%
\begin{picture}(6597,3075)(2596,-5023)
\put(6661,-4021){\makebox(0,0)[lb]{\smash{{\SetFigFont{12}{14.4}{\familydefault}{\mddefault}{\updefault}{\color[rgb]{0,0,0}$\Omega(t)=\plag(\Olag,t)=\pale(\Oale,t)$}%
}}}}
\put(2836,-2986){\makebox(0,0)[lb]{\smash{{\SetFigFont{12}{14.4}{\familydefault}{\mddefault}{\updefault}{\color[rgb]{0,0,0}$\Olag$}%
}}}}
\put(2701,-4876){\makebox(0,0)[lb]{\smash{{\SetFigFont{12}{14.4}{\familydefault}{\mddefault}{\updefault}{\color[rgb]{0,0,0}$\Oale$}%
}}}}
\put(2701,-3976){\makebox(0,0)[lb]{\smash{{\SetFigFont{12}{14.4}{\familydefault}{\mddefault}{\updefault}{\color[rgb]{0,0,0}$\xale$}%
}}}}
\put(2746,-2131){\makebox(0,0)[lb]{\smash{{\SetFigFont{12}{14.4}{\familydefault}{\mddefault}{\updefault}{\color[rgb]{0,0,0}$\xlag$}%
}}}}
\put(5581,-4786){\makebox(0,0)[lb]{\smash{{\SetFigFont{12}{14.4}{\familydefault}{\mddefault}{\updefault}{\color[rgb]{0,0,0}$\xb=\pale \left(\xale,t \right)$}%
}}}}
\put(5401,-2716){\makebox(0,0)[lb]{\smash{{\SetFigFont{12}{14.4}{\familydefault}{\mddefault}{\updefault}{\color[rgb]{0,0,0}$\xb=\plag\left(\xlag,t\right)$}%
}}}}
\put(6661,-3076){\makebox(0,0)[lb]{\smash{{\SetFigFont{12}{14.4}{\familydefault}{\mddefault}{\updefault}{\color[rgb]{0,0,0}$\xb$}%
}}}}
\end{picture}%

%% file: figures/pannello.pdftex_t
\begin{picture}(0,0)%
\includegraphics{figures/pannello.pdf}%
\end{picture}%
\setlength{\unitlength}{2565sp}%
\begingroup\makeatletter\ifx\SetFigFont\undefined%
\gdef\SetFigFont#1#2#3#4#5{%
  \reset@font\fontsize{#1}{#2pt}%
  \fontfamily{#3}\fontseries{#4}\fontshape{#5}%
  \selectfont}%
\fi\endgroup%
\begin{picture}(9849,5574)(514,-8848)
\put(2392,-5757){\makebox(0,0)[lb]{\smash{{\SetFigFont{9}{10.8}{\familydefault}{\mddefault}{\updefault}{\color[rgb]{1,0,0}8}%
}}}}
\put(1186,-8341){\makebox(0,0)[lb]{\smash{{\SetFigFont{9}{10.8}{\familydefault}{\mddefault}{\updefault}{\color[rgb]{0,0,0}1}%
}}}}
\put(3661,-8326){\makebox(0,0)[lb]{\smash{{\SetFigFont{9}{10.8}{\familydefault}{\mddefault}{\updefault}{\color[rgb]{0,0,0}2}%
}}}}
\put(2266,-7036){\makebox(0,0)[lb]{\smash{{\SetFigFont{9}{10.8}{\familydefault}{\mddefault}{\updefault}{\color[rgb]{1,0,0}9}%
}}}}
\put(3710,-7031){\makebox(0,0)[lb]{\smash{{\SetFigFont{9}{10.8}{\familydefault}{\mddefault}{\updefault}{\color[rgb]{1,0,0}7}%
}}}}
\put(3676,-5769){\makebox(0,0)[lb]{\smash{{\SetFigFont{9}{10.8}{\familydefault}{\mddefault}{\updefault}{\color[rgb]{0,0,0}4}%
}}}}
\put(1152,-5764){\makebox(0,0)[lb]{\smash{{\SetFigFont{9}{10.8}{\familydefault}{\mddefault}{\updefault}{\color[rgb]{0,0,0}3}%
}}}}
\put(901,-4861){\makebox(0,0)[lb]{\smash{{\SetFigFont{8}{9.6}{\rmdefault}{\mddefault}{\updefault}{\color[rgb]{0,0,0}$v$}%
}}}}
\put(5026,-8386){\makebox(0,0)[lb]{\smash{{\SetFigFont{8}{9.6}{\rmdefault}{\mddefault}{\updefault}{\color[rgb]{0,0,0}$u$}%
}}}}
\put(8701,-3511){\makebox(0,0)[lb]{\smash{{\SetFigFont{8}{9.6}{\rmdefault}{\mddefault}{\updefault}{\color[rgb]{0,0,0}$ \tb_v^k$}%
}}}}
\put(10126,-5086){\makebox(0,0)[lb]{\smash{{\SetFigFont{8}{9.6}{\rmdefault}{\mddefault}{\updefault}{\color[rgb]{0,0,0}$\tb_u^k$}%
}}}}
\put(3751,-4636){\makebox(0,0)[lb]{\smash{{\SetFigFont{8}{9.6}{\rmdefault}{\mddefault}{\updefault}{\color[rgb]{0,0,0}$\Xb^k(u,v)$}%
}}}}
\put(9676,-3811){\makebox(0,0)[lb]{\smash{{\SetFigFont{8}{9.6}{\rmdefault}{\mddefault}{\updefault}{\color[rgb]{0,0,0}$\nb^k$}%
}}}}
\put(2428,-8333){\makebox(0,0)[lb]{\smash{{\SetFigFont{9}{10.8}{\familydefault}{\mddefault}{\updefault}{\color[rgb]{1,0,0}5}%
}}}}
\put(1134,-7037){\makebox(0,0)[lb]{\smash{{\SetFigFont{9}{10.8}{\familydefault}{\mddefault}{\updefault}{\color[rgb]{1,0,0}6}%
}}}}
\put(9676,-4411){\makebox(0,0)[lb]{\smash{{\SetFigFont{9}{10.8}{\familydefault}{\mddefault}{\updefault}{\color[rgb]{0,0,0}4}%
}}}}
\put(7951,-4036){\makebox(0,0)[lb]{\smash{{\SetFigFont{9}{10.8}{\familydefault}{\mddefault}{\updefault}{\color[rgb]{1,0,0}8}%
}}}}
\put(6601,-5761){\makebox(0,0)[lb]{\smash{{\SetFigFont{9}{10.8}{\familydefault}{\mddefault}{\updefault}{\color[rgb]{1,0,0}6}%
}}}}
\put(6301,-4186){\makebox(0,0)[lb]{\smash{{\SetFigFont{9}{10.8}{\familydefault}{\mddefault}{\updefault}{\color[rgb]{0,0,0}3}%
}}}}
\put(6451,-7186){\makebox(0,0)[lb]{\smash{{\SetFigFont{9}{10.8}{\familydefault}{\mddefault}{\updefault}{\color[rgb]{0,0,0}1}%
}}}}
\put(8026,-7111){\makebox(0,0)[lb]{\smash{{\SetFigFont{9}{10.8}{\familydefault}{\mddefault}{\updefault}{\color[rgb]{1,0,0}5}%
}}}}
\put(9601,-7336){\makebox(0,0)[lb]{\smash{{\SetFigFont{9}{10.8}{\familydefault}{\mddefault}{\updefault}{\color[rgb]{0,0,0}2}%
}}}}
\put(8101,-5611){\makebox(0,0)[lb]{\smash{{\SetFigFont{9}{10.8}{\familydefault}{\mddefault}{\updefault}{\color[rgb]{1,0,0}9}%
}}}}
\put(10051,-5761){\makebox(0,0)[lb]{\smash{{\SetFigFont{9}{10.8}{\familydefault}{\mddefault}{\updefault}{\color[rgb]{1,0,0}7}%
}}}}
\end{picture}%

%% file: figures/pannello.pstex_t
\begin{picture}(0,0)%
\includegraphics{figures/pannello.eps}%
\end{picture}%
\setlength{\unitlength}{2565sp}%
\begingroup\makeatletter\ifx\SetFigFont\undefined%
\gdef\SetFigFont#1#2#3#4#5{%
  \reset@font\fontsize{#1}{#2pt}%
  \fontfamily{#3}\fontseries{#4}\fontshape{#5}%
  \selectfont}%
\fi\endgroup%
\begin{picture}(9849,5574)(514,-8848)
\put(9676,-3811){\makebox(0,0)[lb]{\smash{{\SetFigFont{8}{9.6}{\rmdefault}{\mddefault}{\updefault}{\color[rgb]{0,0,0}$\nb^k$}%
}}}}
\put(9676,-4411){\makebox(0,0)[lb]{\smash{{\SetFigFont{9}{10.8}{\familydefault}{\mddefault}{\updefault}{\color[rgb]{0,0,0}4}%
}}}}
\put(7951,-4036){\makebox(0,0)[lb]{\smash{{\SetFigFont{9}{10.8}{\familydefault}{\mddefault}{\updefault}{\color[rgb]{0,0,0}8}%
}}}}
\put(6601,-5761){\makebox(0,0)[lb]{\smash{{\SetFigFont{9}{10.8}{\familydefault}{\mddefault}{\updefault}{\color[rgb]{0,0,0}6}%
}}}}
\put(6301,-4186){\makebox(0,0)[lb]{\smash{{\SetFigFont{9}{10.8}{\familydefault}{\mddefault}{\updefault}{\color[rgb]{0,0,0}3}%
}}}}
\put(6451,-7186){\makebox(0,0)[lb]{\smash{{\SetFigFont{9}{10.8}{\familydefault}{\mddefault}{\updefault}{\color[rgb]{0,0,0}1}%
}}}}
\put(8026,-7111){\makebox(0,0)[lb]{\smash{{\SetFigFont{9}{10.8}{\familydefault}{\mddefault}{\updefault}{\color[rgb]{0,0,0}5}%
}}}}
\put(9601,-7336){\makebox(0,0)[lb]{\smash{{\SetFigFont{9}{10.8}{\familydefault}{\mddefault}{\updefault}{\color[rgb]{0,0,0}2}%
}}}}
\put(8101,-5611){\makebox(0,0)[lb]{\smash{{\SetFigFont{9}{10.8}{\familydefault}{\mddefault}{\updefault}{\color[rgb]{0,0,0}9}%
}}}}
\put(10051,-5761){\makebox(0,0)[lb]{\smash{{\SetFigFont{9}{10.8}{\familydefault}{\mddefault}{\updefault}{\color[rgb]{0,0,0}7}%
}}}}
\put(2392,-5757){\makebox(0,0)[lb]{\smash{{\SetFigFont{9}{10.8}{\familydefault}{\mddefault}{\updefault}{\color[rgb]{0,0,0}8}%
}}}}
\put(1134,-7037){\makebox(0,0)[lb]{\smash{{\SetFigFont{9}{10.8}{\familydefault}{\mddefault}{\updefault}{\color[rgb]{0,0,0}6}%
}}}}
\put(1186,-8341){\makebox(0,0)[lb]{\smash{{\SetFigFont{9}{10.8}{\familydefault}{\mddefault}{\updefault}{\color[rgb]{0,0,0}1}%
}}}}
\put(2428,-8333){\makebox(0,0)[lb]{\smash{{\SetFigFont{9}{10.8}{\familydefault}{\mddefault}{\updefault}{\color[rgb]{0,0,0}5}%
}}}}
\put(3661,-8326){\makebox(0,0)[lb]{\smash{{\SetFigFont{9}{10.8}{\familydefault}{\mddefault}{\updefault}{\color[rgb]{0,0,0}2}%
}}}}
\put(2266,-7036){\makebox(0,0)[lb]{\smash{{\SetFigFont{9}{10.8}{\familydefault}{\mddefault}{\updefault}{\color[rgb]{0,0,0}9}%
}}}}
\put(3710,-7031){\makebox(0,0)[lb]{\smash{{\SetFigFont{9}{10.8}{\familydefault}{\mddefault}{\updefault}{\color[rgb]{0,0,0}7}%
}}}}
\put(3676,-5769){\makebox(0,0)[lb]{\smash{{\SetFigFont{9}{10.8}{\familydefault}{\mddefault}{\updefault}{\color[rgb]{0,0,0}4}%
}}}}
\put(1152,-5764){\makebox(0,0)[lb]{\smash{{\SetFigFont{9}{10.8}{\familydefault}{\mddefault}{\updefault}{\color[rgb]{0,0,0}3}%
}}}}
\put(901,-4861){\makebox(0,0)[lb]{\smash{{\SetFigFont{8}{9.6}{\rmdefault}{\mddefault}{\updefault}{\color[rgb]{0,0,0}$v$}%
}}}}
\put(5026,-8386){\makebox(0,0)[lb]{\smash{{\SetFigFont{8}{9.6}{\rmdefault}{\mddefault}{\updefault}{\color[rgb]{0,0,0}$u$}%
}}}}
\put(8701,-3511){\makebox(0,0)[lb]{\smash{{\SetFigFont{8}{9.6}{\rmdefault}{\mddefault}{\updefault}{\color[rgb]{0,0,0}$ \tb_v^k$}%
}}}}
\put(10126,-5086){\makebox(0,0)[lb]{\smash{{\SetFigFont{8}{9.6}{\rmdefault}{\mddefault}{\updefault}{\color[rgb]{0,0,0}$\tb_u^k$}%
}}}}
\put(3751,-4636){\makebox(0,0)[lb]{\smash{{\SetFigFont{8}{9.6}{\rmdefault}{\mddefault}{\updefault}{\color[rgb]{0,0,0}$\Xb^k(u,v)$}%
}}}}
\end{picture}%